\documentclass[10pt,a4paper,reqno]{amsart}
\usepackage{clrscode}
\usepackage{mathrsfs,amssymb}

\usepackage{verbatim}

\theoremstyle{plain}

\newtheorem{theorem}{Theorem}

\newtheorem{corollary}{Corollary}
\newtheorem{definition}{Definition}
\newtheorem{proposition}{Proposition}
\newtheorem{lemma}{Lemma}

\theoremstyle{definition}

\newtheorem{remark}{Remark}

\theoremstyle{definition}

\newtheorem{example}{Example}

\newcommand{\im}{\mathrm{im}}

\newcommand{\dotcup}{\ensuremath{\mathaccent\cdot\cup}}

\newcommand{\lb}{\langle}
\newcommand{\rb}{\rangle}

\newcommand{\gr}{\mathscr{R}}
\newcommand{\gl}{\mathscr{L}}
\newcommand{\gh}{\mathscr{H}}

\newcommand{\ZT}{{\mathbb{Z}T}}
\newcommand{\ZS}{{\mathbb{Z}S}}
\newcommand{\ZG}{{\mathbb{Z}G}}
\newcommand{\ZL}{{\mathbb{Z}L}}
\newcommand{\ZH}{{\mathbb{Z}H}}
\newcommand{\ZR}{{\mathbb{Z}R}}
\newcommand{\ZM}{{\mathbb{Z}M}}

\newcommand{\FPn}{{\rm FP}\sb n}
\newcommand{\FPone}{{\rm FP}\sb 1}
\newcommand{\FP}{{\rm FP}}
\newcommand{\FPinfty}{{\rm FP}\sb \infty}

\begin{document}

\title[Homological Finiteness Properties of Monoids]
{Homological finiteness properties of \\ monoids, their ideals and maximal subgroups \\ \medskip \today }

\maketitle

\begin{center}

\vspace{-2mm}

    R. GRAY\footnote{This work was supported by an EPSRC Postdoctoral Fellowship EP/E043194/1 held by the first author at the School of Mathematics \& Statistics of the University of St Andrews, Scotland. \\
    \indent The first author was partially supported by FCT and FEDER, project POCTI-ISFL-1-143 of Centro de \'{A}lgebra da Universidade de Lisboa, and by the project PTDC/MAT/69514/2006.}  
    
\medskip

Centro de \'{A}lgebra da Universidade de Lisboa, \\ Av. Prof. Gama Pinto 2,  1649-003 Lisboa,  Portugal.

    \medskip

    \texttt{rdgray@fc.ul.pt} \\

   \bigskip
   
   \bigskip

    S. J. PRIDE

    \medskip

    Department of Mathematics, University of Glasgow,  \\
    University Gardens, G12 8QW, Scotland.

    \medskip

    \texttt{s.pride@maths.gla.ac.uk} \\

\end{center}

\begin{abstract}
We consider the general question of how the homological finiteness property left-$\FP_n$ (resp. right-$\FPn$) holding in a monoid influences, and conversely depends on, the property holding in the substructures of that monoid. 
This is done by giving methods for constructing free resolutions of substructures from free resolutions of their containing monoids, and vice versa. 
In particular we show that left-$\FP_n$ is inherited by the maximal subgroups in a completely simple minimal ideal, in the case that the minimal ideal has finitely many left ideals. For completely simple semigroups we prove the converse, and as a corollary show that    
a completely simple semigroup is of type left- and right-${\rm FP}\sb n$ if and only if it has finitely many left and right ideals and all of its maximal subgroups are of type ${\rm FP}\sb n$. Also, given an ideal of a monoid, we show that if 
the ideal has a two-sided identity element then the containing monoid is of type left-$\FPn$ if and only if the ideal is of type left-$\FPn$. Applying this result we obtain 
necessary and sufficient conditions for a Clifford monoid (and more generally a strong semilattice of monoids) to be of type left-${\rm FP}\sb n$.
Examples are provided showing that for each of the results all of the hypotheses are necessary.  

\end{abstract}

\section{Introduction}
\label{sec_intro}

Let $S$ be a monoid and $\mathbb{Z}S$ be the monoid ring over the integers $\mathbb{Z}$. For $n \geq 0$ the monoid $S$ is of \emph{type left-${\rm FP}\sb n$} if there is a resolution
\[
A_n \rightarrow A_{n-1} \rightarrow \cdots \rightarrow A_1 \rightarrow A_0 \rightarrow \mathbb{Z} \rightarrow 0
\]
of the trivial left $\mathbb{Z}S$-module $\mathbb{Z}$ such that $A_0, A_1, \ldots, A_n$ are finitely generated free left $\mathbb{Z}S$-modules. 
Monoids of type \emph{right-$\FPn$} are defined dually, working with right $\ZS$-modules. 

The property $\FPn$ was introduced for groups by Bieri in \cite{Bieri1976} and since then has received a great deal of attention in the literature; see \cite{Bestvina1997, Bieri2001, Brady1999, Bux2007, Leary2006}. One natural line of investigation has been the study of the closure properties of $\FPn$. 
Examples include results about the behaviour of $\FPn$ under taking: finite index subgroups or extensions, direct (and semidirect) products, wreath products, HNN extensions, amalgamated free products, and quasi-isometry invariance; see \cite{Alonso1994, Baumslag1998, Bieri1976, Wall1961}. 

In monoid and semigroup theory the property $\FPn$ arises naturally in the 
study of string rewriting systems (i.e. semigroup presentations). 
The history of rewriting systems in monoids and groups is long and distinguished, and has its roots in fundamental work of Dehn and Thue.
The main focus of this research has been on so-called complete rewriting systems (also called convergent rewriting systems) and in algorithms for computing normal forms. A finite complete rewriting system is a finite presentation for a monoid of a particular form (both confluent and Noetherian) which in particular gives a solution of the word problem for the monoid; see \cite{BookAndOtto} for more details. 
Therefore it is of considerable interest to develop an understanding of which monoids are presentable by such rewriting systems.
Many important classes of groups are known to be presentable by finite complete rewriting systems, including Coxeter groups, surface groups, and many closed three-manifold groups.
Rewriting systems continue to receive a lot of attention in the literature; see \cite{Silva2009_2, Silva2009_1, Chouraqui2006, Miasnikov2009, Goodman2008, Hermiller1999, Pride2005}.  The connection between complete rewriting systems and homological finiteness properties is given by a result of
Anick \cite{Anick1986} (see also \cite{Brown1992}) which shows that a monoid that admits such a presentation must be of type left- and right-$\FPinfty$ (meaning type $\FPn$ for all $n$).    
More background on the importance the property $\FP_n$ (and other related finiteness conditions) in semigroup theory, and the connections with the theory of string rewriting systems may be found in the survey articles \cite{Cohen1997, Otto1997}.

For groups the properties left- and right-$\FPn$ are equivalent, so we simply speak of groups of type $\FPn$. However, for monoids in general the two notions are independent. 
Indeed, in \cite{Cohen1992} Cohen gives an example of a monoid, related to the Thompson group, that is right-$\FPinfty$ but not even left-$\FPone$. Several other related homological finiteness properties for monoids have been defined, which are all equivalent to $\FP_n$ when applied to groups, but are different for monoids in general. A central theme of recent research in this area has been to investigate how these various properties relate to one another; see for example \cite{Cohen1992, Cremanns1996, Kobayashi2007, Kobayashi2009, Kobayashi2003, Pride2006, Pride2004}. 

On the other hand, in contrast to the situation in group theory, far less attention has been paid to the closure properties of homological finiteness conditions in semigroup and monoid theory, with only a handful of results of this kind having appeared in the literature. 
For monoid constructions that are direct generalisations of group constructions, perhaps unsurprisingly, some results generalise in a straightforward way from groups to monoids. For example, as observed in \cite{Guba1998}, the direct product $M \times N$ of two monoids is of type left-$\FPn$ if and only if each of $M$ and $N$ is (this may be proved just as for groups using K\"{u}nneth theory). Also generalising from groups, it was shown in \cite{Pride2006} that left-$\FPn$ is inherited under taking retracts of monoids (which is known to be true for groups more generally for quasi-retracts \cite{Alonso1994}). However, the study of constructions specific to, and important in, semigroup theory has not yet received serious attention in the literature.

In recent work \cite{Kobayashi2007, Kobayashi2009} Kobayashi considered the behaviour of left- and right-$\FPn$ for some basic fundamental semigroup-theoretic constructions including left and right zero semigroups (and more generally left and right groups), semilattices, and the process of adjoining a zero element to a monoid. He then used these observations to give examples of monoids that clarify the relationship between the properties left-$\FP_1$, right-$\FP_1$, bi-$\FP_1$, and finite generation (see Section~\ref{sec_FP1} for more on this). The results we obtain here will shed more light on (and in some cases significantly extend) some of the results obtained by Kobayashi in \cite{Kobayashi2007, Kobayashi2009}. 
The importance of understanding the closure properties of $\FPn$ is highlighted further still by the work \cite{Pride2004} where a monoid is constructed from two groups, it is shown how the homological finiteness properties of the monoid relate to those of the groups, and then combined with \cite{Bestvina1997} this is used to give a counterexample to an open question about homological finiteness properties of string rewriting systems.    

In this paper we shall consider the general question of how the property left-$\FP_n$ holding in a monoid influences, and conversely depends on, the property holding in the substructures of that monoid. 
We are particularly interested in relating the property holding in the monoid to the property holding in the subgroups of the monoid, since such results act as a bridge between the homology theory of groups and that of semigroups and monoids.  
The results we present here complement analogous results regarding cohomology obtained in \cite{Adams1967, Nico1969}. 

The paper is structured as follows. 
After giving some preliminaries in Section~\ref{sec_prelims}, we begin our investigation in  
Section~\ref{sec_submonoids} by considering ideals, showing that if $T$ is an ideal of a monoid $S$, and if $T$ has a two-sided identity element, then $S$ is of type left-$\FP_n$ if and only if $T$ is (see Theorem~\ref{thm_IdealWithIdentity}). 
This result is reasonably straightforward to prove. In one direction it generalises the recent observation of Kobayashi \cite{Kobayashi2009} that a monoid with a two-sided zero element is of type left- and right-$\FP_\infty$. 
The result does not extend to ideals in general; see Example~\ref{ex_IdealsInGeneral} in Section~\ref{sec_FP1}. 
Applying this result we show (in Theorem~\ref{thm_Clifford}) that a Clifford monoid $S$, in the sense of \cite{Clifford1941}, is of type left-$\FPn$ if and only if it has a minimal idempotent $e$ and the maximal subgroup that contains $e$ is of type $\FPn$ (for any undefined concepts we refer the reader to Section~\ref{sec_prelims}). 

Results relating properties of monoids with those of their subgroups are known to hold for several other important finiteness conditions including being: finitely generated, finitely presented, having finite derivation type, and being residually finite (see \cite{Golubov1975, GrayMalheiro, Ruskuc1999}).  
In particular for each of these properties it is known that if $S$ is a (von Neumann) regular monoid with finitely many idempotents then $S$ has the property if and only if all of its maximal subgroups have the property.  
In \cite[Remark and Open Problem 4.5]{Ruskuc1999} the author  
posed as an open problem the question of whether the corresponding result holds for the property left-$\FP_n$.
The results mentioned in the previous paragraph answer this question with a resounding no, 
since, for instance, whether a Clifford monoid is of type left-$\FPn$ depends only on one of its maximal subgroups (namely the minimal one), and the other maximal subgroups can have any properties that one desires. 
This leaves the general question of to what extent the property left-$\FP_n$ holding in a monoid relates to the property holding in the maximal subgroups of that monoid. As we shall see, it is straightforward to show (see Theorem~\ref{thm_MinimalMaximalSubgroups}) that the property left-$\FP_n$ is inherited by maximal subgroups contained in a completely simple minimal ideal, in the case that the ideal has finitely many left ideals. Without the finiteness assumption on left ideals the result no longer holds (see Example~\ref{ex_NeedBFinite} in Section~\ref{sec_applications}).  
In Sections~\ref{sec_resolutions} and \ref{sec_resolutions2} we consider the more difficult converse problem, and present results that show how the property left-$\FP_n$ holding in a completely simple semigroup relates to 
the property holding in its maximal subgroups (see Theorem~\ref{thm_completelysimplemain}, Corollary~\ref{corollary_main1} and Theorem~\ref{thm_main2}). In particular we show that  
a completely simple semigroup is of type left- and right-${\rm FP}\sb n$ if and only if it has finitely many left and right ideals and all of its maximal subgroups are of type ${\rm FP}\sb n$.

In Section~\ref{sec_FP1} using recent results of Kobayashi \cite{Kobayashi2007} we go on to analyse left-$\FP_1$ in more detail and, extending \cite[Corollary~2.7]{Kobayashi2007}, we give necessary and sufficient conditions for a completely simple semigroup to be of type left-$\FP_1$ (see Theorem~\ref{thm_fp1compsimple}).  
As an application we deduce that a monoid with finitely many left and right ideals is of type left-$\FP_1$ if and only if the maximal subgroups contained in its (unique) minimal ideal are all finitely generated. 
In Section~\ref{sec_applications} we give examples showing the hypotheses of our main results are necessary, and present some further applications. Finally, in Section~\ref{sec_other_props} we discuss some other related homological finiteness properties, including having finite cohomological dimension, and explain how one may construct counterexamples to several other open problems posed in \cite{Ruskuc1998} and \cite{Ruskuc1999} regarding these properties.   
 
\section{Preliminaries}
\label{sec_prelims}

\subsection*{\boldmath Free resolutions and the finiteness property $\FP_n$}

Let us begin by recalling some basic definitions from homology theory that we need; for more details we refer the reader to \cite{Hilton1997}.

Let $S$ be a monoid and let $\ZS$ be the integral monoid ring of $S$. 
We have the standard augmentation
\[
\epsilon_S: \ZS \rightarrow \mathbb{Z}, \quad s \mapsto 1 \quad (s \in S)
\]
and therefore we can regard $\mathbb{Z}$ as a trivial left $\ZS$-module with the $\ZS$-action via $\epsilon_S$:
\[
\lambda \cdot z = \epsilon_S(\lambda) z \quad (\lambda \in \ZS, \; z \in \mathbb{Z}).
\]
A \emph{free resolution} of the trivial left $\ZS$-module $\mathbb{Z}$ is a sequence $A_0, A_1, A_2, \ldots$ of free left $\ZS$-modules and homomorphisms $\partial_0: A_0 \rightarrow \mathbb{Z}$ and $\partial_i: A_i \rightarrow A_{i-1}$, for $n \geq 1$, such that the sequence
\[
\cdots \rightarrow
A_2 \xrightarrow{\partial_2}
A_1 \xrightarrow{\partial_1}
A_0 \xrightarrow{\partial_0}
\mathbb{Z} \rightarrow 0
\]
is \emph{exact} (i.e. $\im \partial_{n+1} = \mathrm{ker} \partial_n$ for $n \geq 0$, and $\partial_0(A_0) = \mathbb{Z}$).
We shall often refer to such a resolution simply as a \emph{free left resolution of $S$}.
A monoid
$S$ is said to be of type left-$\FP_n$ if there is a partial free resolution
of the trivial left $\ZS$-module $\mathbb{Z}$:
\[
A_n \rightarrow A_{n-1} \rightarrow \cdots \rightarrow A_1 \rightarrow A_0 \rightarrow \mathbb{Z} \rightarrow 0
\]
where $A_0, A_1, \ldots, A_n$ are all finitely generated. Dually we can regard $\mathbb{Z}$ as a right $\ZS$-module via $\epsilon_S$, and analogously define free right resolutions, and monoids of type right-$\FP_n$.
For a semigroup $S$ without a two-sided identity the ring $\ZS$ does not have an identity and the above definition is no longer valid. Thus to extend the notion of $\FPn$ to arbitrary semigroups we utilise the standard device of adjoining an identity element. 
That is, given a semigroup $S$ we use $S^1$ to denote the semigroup $S = S \cup \{ 1 \}$ with an identity element $1$ adjoined with $1 \not\in S$ and where we define $s1 = 1s = s$ for all $s \in S$.
Then for a semigroup $S$ without a two-sided identity element we say that $S$ is of type left-$\FPn$ if the monoid $S^1$ is of type left-$\FPn$.

There are several important differences when working with $\ZS$-modules, where $S$ is a monoid, compared to working with $\ZG$-modules, with $G$ a group. Two of the most important differences are the following. 

\vspace{2mm}

\noindent $\bullet$  For a group $G$, any left $\ZG$-module can be regarded as a right $\ZG$-module, and conversely, by defining $u g = g^{-1} u$ for any $g \in G$ and $u$ in the module. Moreover applying this operation to a free left resolution of $\mathbb{Z}$ yields a free right resolution of $\mathbb{Z}$ and hence (as mentioned in the introduction)
for groups left-$\FP_n$ and right-$\FP_n$ are equivalent, and we simply speak of property $\FP_n$ when working with groups. More generally left- and right-$\FP_n$ coincide for inverse semigroups (see \cite{Lawson1998} for more on inverse semigroups).
However, as observed in the introduction, the same is far from being true for arbitrary monoids; see \cite{Cohen1992}. 

\vspace{2mm}

\noindent $\bullet$ If $G$ is a group and $H$ is a subgroup of $G$ then $\ZG$, regarded as a left module over $\ZH$, is free with rank equal to the index of $H$ in $G$.
This fact is fundamental for the proof that $\FPn$ is preserved under taking finite index subgroups or extensions (see \cite[Proposition~5.1]{Brown1982}).
In contrast, given a monoid $S$ and subsemigroup $T$, where $T$ has with a two-sided identity element, then $\ZS$ regarded as a left module over $\ZT$ will \emph{not} in general be free. 

\vspace{2mm}

\noindent These two fundamental differences are the main reason for the fairly intricate arguments needed to establish some of the results below.

We shall make repeated use of the following consequence of the generalised Schanuel Lemma; see {\cite[p193]{Brown1982}}.
\begin{lemma}
\label{lem_Schanuel}
Let $S$ be a monoid. For $n \geq 0$ if
\[
A_{n-1} \xrightarrow{\partial_{n-1}}
\cdots
\xrightarrow{\partial_{2}}
A_1
\xrightarrow{\partial_{1}}
A_0
\xrightarrow{\partial_{0}}
\mathbb{Z}
\rightarrow 0
\]
is a partial free resolution of the trivial left $\ZS$-module $\mathbb{Z}$, with 
$A_0, \ldots, A_{n-1}$ finitely generated free left $\ZS$-modules, then $S$ is of type left-${\rm FP}\sb n$ if and only if $\ker \partial_{n-1}$ is finitely generated.
\end{lemma}
One consequence of Lemma~\ref{lem_Schanuel} is that if a monoid $S$ is of type left-$\FP_n$ then there is a partial free resolution:
\[
A_{n-1} \xrightarrow{\partial_{n-1}}
\ldots
\xrightarrow{\partial_{2}}
A_1
\xrightarrow{\partial_{1}}
A_0
\xrightarrow{\partial_{0}}
\mathbb{Z}
\rightarrow 0
\]
of the trivial left $\ZS$-module $\mathbb{Z}$ where $A_0 = \ZS$ and $\partial_0 = \epsilon_S$ is the standard augmentation. 

Throughout, given a subset $X$ of a left $\ZS$-module $A$, where $S$ is a monoid, we use $\lb X \rb_\ZS$ to denote the left $\ZS$-module generated by the set $X$. Unless otherwise stated, we work with left modules and property left-$\FP_n$ throughout. 

\subsection*{Green's relations and completely simple semigroups}

We now outline some of the basic concepts from semigroup theory that we shall need; for more details we refer the reader to \cite{Howie1995, Rhodes2009}. 
Green's relations were first introduced in \cite{Green1951} and have ever since played a fundamental role in the structure theory of semigroups. For elements $x$ and $y$ of a semigroup $S$ we write $x \gr y$ if $x$ and $y$ generate the same principal right ideal, $x \gl y$ if they generate the same principal left ideal, and let $\gh$ denote the intersection of $\gr$ and $\gl$. In other words, for $x,y \in S$
\[
x \gr y \Leftrightarrow xS^1 = yS^1, \quad
x \gl y \Leftrightarrow S^1 x = S^1 y, \quad
x \gh y \Leftrightarrow x \gl y \wedge x \gr y.
\]
Each of these relations is an equivalence relation on $S$ which we call the $\gr$-, $\gl$- and $\gh$-classes of the semigroup, respectively. A semigroup $S$ is said to be (von Neumann) \emph{regular} if for every $x\in S$ there exists $y\in S$ such that $xyx=x$.
A semigroup is regular if and only if every $\gr$-class (equivalently every $\gl$-class) contains at least one idempotent.

We use $E(S)$ to denote the set of idempotents of a semigroup $S$.
Let $e$ be an idempotent in a semigroup $S$. Then $eSe$ is the largest submonoid of $S$ (with respect to inclusion) whose identity element is $e$. The group of units $G_e$ of $eSe$ (i.e. the members of $eSe$ that have two-sided inverses in $eSe$) is the largest subgroup of $S$ (with respect to inclusion) with identity $e$, and is called the \emph{maximal subgroup of $S$ containing $e$}. If an $\gh$-class $H$ contains an idempotent then $H$ is a maximal subgroup of $S$, and conversely every maximal subgroup of $S$ arises in this way.
From the definitions it is easily seen that the $\gr$-classes (resp. $\gl$-classes) are in one-one correspondence with principal right (resp. left) ideals of the semigroup.

Of particular importance are those semigroups that have no proper two-sided ideals. 
A semigroup is called \emph{simple} if it has no proper two-sided ideals, and is called \emph{completely simple} if it is simple and has minimal left and right ideals.
It is not hard to see that if a semigroup $S$ has a minimal ideal $K$ then $K$ is a simple semigroup, which is sometimes referred to as the \emph{kernel} of the semigroup. 
A \emph{right zero semigroup} is a semigroup $U$ such that $xy = y$ for all $x,y \in U$ (dually one defines left zero semigroup), and a right (resp. left) group is a direct product of a group and a right (resp. left) zero semigroup. A right (resp. left) group is precisely a completely simple semigroup with a single $\gr$-class (resp. $\gl$-class).
A \emph{rectangular band} is a direct product of a left zero semigroup and a right zero semigroup.  
Free resolutions for completely simple semigroups will be considered in detail in 
Sections~\ref{sec_resolutions} and \ref{sec_resolutions2}. 

\section{Resolutions for Subsemigroups}
\label{sec_submonoids}

In this section we make some general observations about the relationship between the property left-$\FP_n$ holding in a monoid and the same property holding in certain subsemigroups of the monoid. 
We shall see that if $M$ is an ideal of a monoid $S$, and $M$ has a two-sided identity element, then passing from free resolutions of $S$ to free resolutions of $M$ is straightforward. In particular, we shall prove the following result.

\begin{theorem}
\label{thm_McrossB}
Let $S$ be a monoid, let $R$ be a right ideal of $S$ and suppose that $R \cong M \times B$ where $M$ is a monoid and $B$ is a right zero semigroup. If $S$ is of type left-$\FP_n$ and $B$ is finite then $M$ is of type left-$\FP_n$. 
\end{theorem}

Once established, Theorem~\ref{thm_McrossB} can then be applied both to maximal subgroups in minimal ideals (Theorem~\ref{thm_MinimalMaximalSubgroups}) and to ideals with identity (Theorem~\ref{thm_IdealWithIdentity}).    
Before proving Theorem~\ref{thm_McrossB} we first need some basic lemmas. 

Let $S$ be a monoid and let $M$ be a subsemigroup of $S$ such that $M$ has a two-sided identity element $e \in M$. If $A$ is a left $\ZS$-module then $e \cdot A = eA$ is a left $\ZM$-module with
\[
ea_1 + ea_2 = e(a_1 + a_2) \ (a_1, a_2 \in A), 
\
\lambda (e a) = e (\lambda a) \in eA 
\
(a \in A, \lambda \in \ZM). 
\]
There is then an obvious functor $\Phi$ from the category of left $\ZS$-modules to the category of left $\ZM$-modules defined as follows. For a left $\ZS$-module $A$ we define $\Phi(A) = eA$, and for a left $\ZS$-module homomorphism $\theta: A_2 \rightarrow A_1$ we let
\[
\Phi(\theta): eA_2 \rightarrow eA_1
\]
be the restriction of $\theta$ to $eA_2$. This is well defined since for all $a_2 \in A_2$ we have
$
\theta(ea_2) = e \theta(a_2) \in eA_1. 
$
It is not hard to see that the functor $\Phi$ is exact, so we omit the proof. 
\begin{lemma}
\label{lem_exactness}
The functor $\Phi$ is exact i.e. if
\[
A_2 \xrightarrow{\theta_2} A_1 \xrightarrow{\theta_1} A_0
\]
is an exact sequence of left $\ZS$-modules then
\[
\Phi(A_2) \xrightarrow{\Phi(\theta_2)} \Phi(A_1) \xrightarrow{\Phi(\theta_1)} \Phi(A_0)
\]
is an exact sequence of left $\ZM$-modules. 
\end{lemma}
For a general subsemigroup $M$, with a two-sided identity element, of a semigroup $S$ the functor $\Phi$ will not map free left $\ZS$-modules to free left $\ZM$-modules. We now show that when $S$ and $M$ satisfy the conditions given in the statement of Theorem~\ref{thm_McrossB}, then freeness is preserved.

\begin{lemma}
\label{lem_free_inS}
Let $R = N \times B$ where $N$ is a monoid and $B$ is a right zero semigroup. 
Fix $y \in B$ and let $M = \{(n,y) : n \in N  \}$.  
Then, viewed as a left $\ZM$-module, $\ZR$ is free with basis $F = \{ (1,b) : b \in B \}$. 
\end{lemma}
\begin{proof}
Clearly each $r \in R$ can be written uniquely in the form $r=mf$ where $m \in M$ and $f \in F$.  
It follows that each $\alpha \in \ZR$ can be written uniquely in the form $\alpha = \sum_{f \in F} \lambda_f f$ where $\lambda_f \in \ZM$ for $f \in F$. 
This proves the lemma. 
\end{proof}

\begin{lemma}
\label{lem_freeness}
Let $S$ be a monoid with a right ideal $R$ such that $R = N \times B$ where $N$ is a monoid and $B$ is a right zero semigroup. Fix $y \in B$ and let $M = \{ (n,y): n \in N \}$. Then $M$ is a subsemigroup of $S$ with a two-sided identity $e = (1,y)$, and if $A$ is a free left $\ZS$-module of rank $r$ then $\Phi(A) = eA$ is a free left $\ZM$-module of rank $r|B|$.
\end{lemma}
\begin{proof}
Let $A = \bigoplus_{x \in X} \ZS x$ be a free left $\ZS$-module with basis $X$, where $X$ is a non-empty set with $|X|=r$. It follows from the hypotheses that $eS = R$ and therefore $eA = \bigoplus_{x \in X} \ZR x$ which is a free left $\ZM$-module with basis
\[
F \cdot X = \{ f \cdot x : f \in F, x \in X  \}
\]
by Lemma~\ref{lem_free_inS}, where $F = \{ (1,b) : b \in B \}$. Therefore $eA$ is a free left $\ZM$-module of rank $|X||B| = r|B|$.    
\end{proof} 

\begin{proof}[Proof of Theorem~\ref{thm_McrossB}]
Suppose $S$ is of type left-$\FPn$ and that $B$ is finite. Let 
\[
\mathcal{A}:
A_n \xrightarrow{\theta_n}
A_{n-1} \xrightarrow{\theta_{n-1}}
\cdots \xrightarrow{\theta_2} 
A_1 \xrightarrow{\theta_1} 
A_0 \xrightarrow{\theta_0} 
\mathbb{Z} \rightarrow 0
\]
be a partial free left resolution for $S$ where $A_i$ is a finitely generated free left $\ZS$-module for $i = 0, \ldots, n$. It then 
follows from Lemmas~\ref{lem_exactness} and \ref{lem_freeness} that
\[
\mathcal{B}:
B_n \xrightarrow{\psi_n}
B_{n-1} \xrightarrow{\psi_{n-1}}
\cdots \xrightarrow{\psi_2} 
B_1 \xrightarrow{\psi_1} 
B_0 \xrightarrow{\psi_0} 
\mathbb{Z} \rightarrow 0
\]
is a partial free left resolution for $M$ where $B_i = eA_i$ is a finitely generated left $\ZM$-module
and $\psi_i$ is the restriction of $A_i$ to $B_i$, for $i=0, \ldots, n$. Therefore $M$ is of type left-$\FPn$.
\end{proof}

We can apply Theorem~\ref{thm_McrossB} to obtain the following. 

\begin{theorem}
\label{thm_MinimalMaximalSubgroups}
Let $S$ be a monoid and let $H$ be a maximal subgroup of $S$ contained in a completely simple minimal ideal $U$ of $S$. If $S$ is of type left-$\FP_n$ and $U$ has finitely many left ideals then $H$ is of type $\FP_n$.
\end{theorem}
\begin{proof}
Let $e \in H$ be the identity of $H$, let $R = eS$ and set $F = E(S) \cap R$. 
Since $U$ is a minimal ideal and is completely simple 
it follows that $R$ is an $\gr$-class of $U$ which by the Rees theorem \cite[Section~3.2]{Howie1995} implies that $R \cong H \times F$. But $H$ is a monoid and $F$, which is a set of $\gr$-related idempotents, is a right zero semigroup. Now the result follows by applying Theorem~\ref{thm_McrossB}.   
\end{proof}

The converse of Theorem~\ref{thm_MinimalMaximalSubgroups} does not hold. Indeed, if $L$ is an infinite left zero semigroup, then $S=L^1$ has finitely many left ideals and all of its maximal subgroups are trivial (and so are of type left- and right-$\FP_\infty$) but $S$ itself is not of type left-$\FP_1$ by Theorem~\ref{thm_fp1compsimple} below. The same example shows that the converse of Theorem~\ref{thm_McrossB} is also not true in general. 

Theorem~\ref{thm_MinimalMaximalSubgroups} may, in particular, be applied to completely simple semigroups with finitely many left ideals. Necessary and sufficient conditions for such a semigroup to be of type left-$\FP_n$ will be given in Theorem~\ref{thm_completelysimplemain} below. 

The assumption that $B$ is finite is necessary for Theorem~\ref{thm_McrossB} to hold (correspondingly the assumption that $U$ has finitely many left ideals is necessary for Theorem~\ref{thm_MinimalMaximalSubgroups}); see Example~\ref{ex_NeedBFinite} in Section~\ref{sec_applications} below. When $B$ is a singleton, $M$ is an ideal with a two-sided identity, and in this case the converse of Theorem~\ref{thm_McrossB} does hold, as we now demonstrate.

\begin{theorem}
\label{thm_IdealWithIdentity}
Let $S$ be a monoid, let $T$ be an ideal of $S$ and suppose that $T$ has a two-sided identity element. 
Then $S$ is of type left-$\FP_n$ if and only if $T$ is of type left-$\FP_n$.
\end{theorem}
\begin{proof}
Applying Theorem~\ref{thm_McrossB} in the case $|B|=1$ proves that if $S$ is of type left-$\FP_n$ then so is $T$. (Alternatively, this direction follows from \cite[Theorem~3]{Pride2006}, since $T$ is a retract of $S$.)

For the converse,
suppose that $T$ is of type left-$\FP_n$.
By Lemma~\ref{lem_Schanuel} this means that there is
a partial free resolution
\[
\mathcal{A}:
A_n \xrightarrow{\partial_n}
A_{n-1} \xrightarrow{\partial_{n-1}}
\cdots \xrightarrow{\partial_2} 
A_1 \xrightarrow{\partial_1} 
A_0 \xrightarrow{\partial_0 = \epsilon_T} \mathbb{Z} \rightarrow 0
\]
of the trivial left $\ZT$-module $\mathbb{Z}$ where $A_0 = \ZT$, 
$\partial_0 = \epsilon_T$ is the standard augmentation,
$\ker \partial_i = \lb X_i \rb_\ZT$ with $X_i \subseteq A_i$ and $X_i$ finite for all $0 \leq i \leq n-1$
where 
\[
A_i = \bigoplus_{x \in X_{i-1}} \ZT [x],
\]
for $1 \leq i \leq n$, 
and $\partial_i:A_i \rightarrow A_{i-1}$ is the left $\ZT$-module homomorphism extending $[x] \mapsto x$ for $x \in X_{i-1}$. 
Using $\mathcal{A}$, we now construct a partial free resolution for $S$. 
Let $e \in T$ be the two-sided identity of $T$.
Define:
\[
\mathcal{B}:
B_n \xrightarrow{\partial_n'}
B_{n-1} \xrightarrow{\partial_{n-1}'}
\cdots \xrightarrow{\partial_2'}
B_1
\xrightarrow{\partial_1'}
B_0 \xrightarrow{\partial_0' = \epsilon_S} \mathbb{Z} \rightarrow 0
\]
where $B_0 = \ZS$ and for $1 \leq i \leq n$:
\[
B_i =
\bigoplus_{j=i-1, \ldots, 0 \atop x \in X_j} \ZS[x]
\oplus \ZS [e], 
\]
and where, viewing $A_i$ as a subset of $B_i$ under the obvious inclusion arising from $\ZT \subseteq \ZS$, the map $\partial_i' : B_i \rightarrow B_{i-1}$ is given by
\begin{eqnarray*}
\partial_i'([x]) & = &
\begin{cases}
x & \mbox{if $x \in X_{i-1}$} \\
(1-e)[x] & \mbox{if $x \in X_j$ $(j = i-2, i-4, \ldots)$ } \\
e[x] & \mbox{if $x \in X_j$ $(j=i-3, i-5, \ldots)$,}
\end{cases} \\
\partial_i'([e]) & = &
\begin{cases}
e[e] & \mbox{if $i$ is even} \\
(1-e) & \mbox{if $i=1$} \\
(1-e)[e] & \mbox{if $i \neq 1$ and $i$ is odd.} 
\end{cases}
\end{eqnarray*}

\smallskip

\noindent \textbf{Claim~1.} Let $Y_0 = X_0 \cup \{ (1-e) \}$ and for all $1 \leq k \leq n-1$ let
\begin{eqnarray*}
Y_k = 
X_k & \cup &
\{
(1-e)[x]: x \in X_i, \; i = k-1, k-3, \ldots
\} \\
& \cup &
\{
e[x]: x \in X_i, \; i = k-2, k-4, \ldots
\}
\cup Z
\end{eqnarray*}
where
\[
Z =
\begin{cases}
e[e] & \mbox{if $k$ is odd} \\
(1-e)[e] & \mbox{if $k$ is even}.
\end{cases}
\]
Then $Y_k$ is a subset of $\ker \partial_k'$ and $\lb Y_k \rb_\ZS = \ker \partial_k'$, for all $0 \leq k \leq n-1$. 
\begin{proof}[Proof of Claim~1.]
We shall consider only the case that $k$ is odd. The case $k$ even, and in particular the case $k=0$, may be dealt with using a similar argument. 

To see that $Y_k \subseteq \ker \partial_k'$ first note that $X_k \subseteq \ker \partial_k \subseteq \ker \partial_k'$ under the natural inclusion $A_k \subseteq B_k$ arising from $\ZT \subseteq \ZS$.
Also
\[
\partial_k'((1-e)[x]) = (1-e)x = x-x = 0
\]
for all $x \in X_{k-1}$, since 
\begin{equation}
\label{eqn_star_star}
X_{k-1} \subseteq \ker \partial_{k-1} \subseteq A_{k-1} = \bigoplus_{x \in X_{k-2}} \ZT[x]
\end{equation}
and $e$ is a two-sided identity for $T$ (when $k=1$, $X_{k-1} = X_0 \subseteq \ZT$ and we also have $(1-e)x=0$ for all $x \in X_{k-1}$).  
It is then easily checked that the remaining members of $Y_k$ all belong to $\ker \partial_k'$.

We are left with the task of proving $\lb Y_k \rb_\ZS = \ker \partial_k'$. To this end, let $\alpha \in \ker \partial_k'$ be arbitrary, say
\[
\alpha =
\sum_{i = k-1, \ldots, 0 \atop x \in X_i}
\lambda_x [x]
+ \lambda_e [e]
\]
where $\lambda_e \in \ZS$ and each $\lambda_x \in \ZS$.
Since $\alpha \in \ker \partial_k'$ and $k$ is odd we have
\[
0 
=
\partial_k'(\alpha)
=
\sum_{x \in X_{k-1}}
\lambda_x x +
\sum_{i=k-3, k-5, \ldots \atop x \in X_i}
\lambda_x e [x] +
\sum_{i=k-2, k-4, \ldots \atop x \in X_i}
\lambda_x (1-e) [x] +
\lambda_e (1-e)[e].
\]
(When $k=1$ the last term in the above expression will simply be $\lambda_e(1-e)$.)
Along with \eqref{eqn_star_star}, consideration of the coefficients of this equation gives
\[
\begin{array}{rl}
\lambda_x e = 0 			&  \mbox{for} \ x \in X_i, \; i=k-3, k-5, \ldots \\
\lambda_x (1-e) = 0 	&  \mbox{for} \ x \in X_i, \; i=k-4, k-6, \ldots \\
\lambda_e (1-e) = 0. &
\end{array}
\]
Therefore:
\begin{eqnarray*} 
\sum_{i=k-3, k-4, \ldots \atop x \in X_i}
\lambda_x [x] + \lambda_e [e] 
& = & 
\sum_{i=k-3, k-5, \ldots \atop x \in X_i}
\lambda_x [x]  
+
\sum_{i=k-4, k-6, \ldots \atop x \in X_i}
\lambda_x [x]  
+
\lambda_e [e] \\
& = &
\sum_{i=k-3, k-5, \ldots \atop x \in X_i}
\lambda_x (1-e) [x]  
+
\sum_{i=k-4, k-6, \ldots \atop x \in X_i}
\lambda_x e [x]  
+
\lambda_e e [e] \\
&\in&
\lb Y_k \rb_\ZS,
\end{eqnarray*}
by inspection of $Y_k$.  
In other words $\alpha - \alpha_1 \in \lb Y_k \rb_\ZS$ where
\[
\alpha_1 =
\sum_{x \in X_{k-1}}
\lambda_x [x]
+
\sum_{x \in X_{k-2}}
\lambda_x [x]
\in \ker \partial_k',
\]
and
\begin{equation} 
\label{eqn_alpha1}
0 = \partial_k'(\alpha_1)  = 
\sum_{x \in X_{k-1}} \lambda_x x
+
\sum_{x \in X_{k-2}} \lambda_x (1-e) [x].
\end{equation} 
(When $k=1$ equation \eqref{eqn_alpha1} will actually be $0 = \sum_{x \in X_0} \lambda_x x + \lambda_e (1-e)$, which implies $\lambda_e \in \ZT$, and then the rest of the argument follows the same lines as below.)
Since $T$ is an ideal and $X_{k-1} \subseteq A_{k-1}$, it follows that 
$\sum_{x \in X_{k-1}} \lambda_x x \in A_{k-1}$ which by equation \eqref{eqn_alpha1} implies  
$\lambda_x (1-e) \in \ZT$ for all $x \in X_{k-2}$. But clearly this is only possible if $\lambda_x \in \ZT$ for all $x \in X_{k-2}$. Therefore:
\[
\sum_{x \in X_{k-2}} \lambda_x [x] 
=
\sum_{x \in X_{k-2}} \lambda_x e [x]
\in
\lb Y_k \rb_\ZS,
\]
and this implies $\alpha_1 - \alpha_2 \in \lb Y_k \rb_\ZS$ where
\[
\alpha_2 = \sum_{x \in X_{k-1}} \lambda_x [x] \in \ker \partial_k'.
\]
Now
\[
\sum_{x \in X_{k-1}} \lambda_x [x] -
\sum_{x \in X_{k-1}} \lambda_x e [x] 
=
\sum_{x \in X_{k-1}} \lambda_x (1-e) [x] 
\in \lb Y_k \rb_\ZS.
\]
But since $\alpha_2 \in \ker \partial_k'$ and $e$ is a left identity for $T$ it follows that
\[
\sum_{x \in X_{k-1}} \lambda_x e [x] \in \ker \partial_k = \lb X_k \rb_\ZT \subseteq \lb Y_k \rb_\ZS,
\]
and thus $\alpha_2 \in \lb Y_k \rb_\ZS$. Combined with the previous observations above we conclude
that 
$
\alpha = (\alpha - \alpha_1) + (\alpha_1 - \alpha_2) + \alpha_2
$
belongs to $\lb Y_k \rb_\ZS$, completing the proof of the claim.
\end{proof}
Returning to the proof of Theorem~\ref{thm_IdealWithIdentity}, by definition each $B_k$ is a finitely generated free left $\ZS$-module and for $1 \leq k \leq n$ the mapping $\partial_k':B_k \rightarrow B_{k-1}$ is a left $\ZS$-module homomorphism. 
From the definitions it is easily seen that $Y_k$ is a subset of $\im \partial_{k+1}'$ 
and that $\im \partial_{k+1}' \subseteq \lb Y_k \rb_\ZS$. Combined with Claim~1, this shows that $\mathcal{B}$ is exact, completing the proof of the theorem. 
\end{proof}

Kobayashi \cite{Kobayashi2009} recently observed that a monoid with a two-sided zero element is of type left- and right-$\FP_\infty$. Kobayashi's result is a special case of Theorem~\ref{thm_IdealWithIdentity} where $T = \{ 0 \}$. 

An analogous result to Theorem~\ref{thm_IdealWithIdentity} regarding cohomology was proved in \cite{Adams1967}.

Example~\ref{ex_IdealsInGeneral} in Section~\ref{sec_FP1} below shows that Theorem~\ref{thm_IdealWithIdentity} does not hold if we remove the assumption that the ideal $T$ has a two-sided identity element. 
Theorem~\ref{thm_IdealWithIdentity} will be applied below in Section~\ref{sec_applications} to prove Theorem~\ref{thm_Clifford} which characterises the property left-$\FPn$ for Clifford monoids (and more generally strong semilattices of monoids). 

\section{Resolutions for Completely Simple Semigroups I}
\label{sec_resolutions}

Let $U$ be a completely simple semigroup, let $H$ be a maximal subgroup of $U$, and let $S = U^1$. 
As we saw above in Section~\ref{sec_submonoids}, 
given a free left resolution for $S$ we may construct a free left resolution for $H$. In particular, if the free left resolution for $S$ is finitely generated up to dimension $n$, and $U$ has only finitely many left ideals, then the free left resolution for $H$ will be finitely generated up to dimension $n$ also i.e. the property left-$\FP_n$ will be inherited by $H$ from $S$. 

In this section and the one that follows it we shall consider the converse problem. Given a partial free resolution for the group $H$ we show how to construct a partial free resolution for $S=U^1$, and then use this to prove the following. 

\begin{theorem}
\label{thm_completelysimplemain}
Let $U$ be a completely simple semigroup with finitely many left ideals and let $H$ be a maximal subgroup of $U$. Then $U$ is of type left-${\rm FP}\sb n$ if and only if $U$ has finitely many right ideals and the group $H$ is of type ${\rm FP}\sb n$.
\end{theorem}

Examples will be provided in Sections \ref{sec_FP1} and ~\ref{sec_applications} showing that this theorem fails if any of the hypotheses are lifted. 
Theorem~\ref{thm_completelysimplemain} has the following immediate consequence.

\begin{corollary}
\label{corollary_main1}
Let $U$ be a finitely generated completely simple semigroup. Then the following are equivalent:
\begin{enumerate}
\item[(i)] $U$ is of type left-$\FP_n$;
\item[(ii)] $U$ is of type right-$\FP_n$;
\item[(iii)] All maximal subgroups of $U$ are of type $\FP_n$.
\end{enumerate}
\end{corollary}

We also have the following. 

\begin{sloppypar}
\begin{theorem}
\label{thm_main2}
Let $U$ be a completely simple semigroup and let $H$ be a maximal subgroup of $U$.
Then $U$ is of type left-$\FP_n$ and right-${\rm FP}\sb n$ if and only if $U$ has finitely many left and right ideals and the group $H$ is of type ${\rm FP}\sb n$.
\end{theorem}
\end{sloppypar}
\begin{proof}
Suppose that $U$ is of type left- and right-$\FPn$. Then in particular $U$ is of type left- and right-$\FP_1$ which, by Theorem~\ref{thm_fp1compsimple} (or alternatively \cite[Theorem~2.6]{Kobayashi2007}), implies that $U$ has finitely many left and right ideals. Then by Theorem~\ref{thm_completelysimplemain} it follows that $H$ is of type $\FPn$. 

The converse is a direct corollary of Theorem~\ref{thm_completelysimplemain}.  
\end{proof}

The rest of this section, and the one that follows it, will be dedicated to the proof of Theorem~\ref{thm_completelysimplemain}. 
Results analogous to those above for other finiteness properties, including automaticity and finite derivation type, have appeared in the literature; see \cite{Campbell2002, Descalco2002, Malheiro2006}.
Other important recent work on completely simple semigroups includes \cite{Kambites2007_2}. 

Let us outline our strategy for proving Theorem~\ref{thm_completelysimplemain}.  
We shall adopt standard notation for completely simple semigroups.
Let $U$ be a completely simple semigroup. 
We assume that the $\gr$- and $\gl$-classes of $U$ are indexed by sets $I$ and $\Omega$ respectively so that
\[
U = \bigcup_{i \in I} R_i = \bigcup_{\omega \in \Omega} L_{\omega}.
\]
The $\gh$-classes of $U$ are the sets $H_{i \omega} = R_i \cap L_{\omega}$ for $i \in I$ and $\omega \in \Omega$. Every $\gh$-class of $U$ contains an idempotent, we use $e_{i \omega}$ to denote the idempotent of $H_{i \omega}$ which is exactly the identity of the group $H_{i \omega}$. All the group $\gh$-classes $H_{i \omega} \;(i \in I, \omega \in \Omega)$ are isomorphic to a fixed group $G$, called the \emph{Sch\"{u}tzenberger group} of $U$. 
The best way to visualise a completely simple semigroup $U$ is as a rectangular grid tiled with $|I| \times |\Omega|$ squares, representing the $\gh$-classes, with each row of squares representing an 
$\gr$-class, and each column of squares representing an $\gl$-class 
(the is sometimes referred to as an \emph{egg box diagram}). 

The following proposition lists some basic properties of completely simple semigroups.

\begin{proposition}
\label{prop_compsimpleprops}
Let $U$ be a completely simple semigroup with set of $\gr$-, $\gl$- and $\gh$-classes 
$
\{ R_i : i \in I \}$,
$\{ L_\omega : \omega \in \Omega \}$ and
$\{ H_{i \omega} : i \in I, \omega \in \Omega \}
$, respectively. 
\begin{enumerate}
\item[(i)] If $x \in H_{i \omega}$ and $y \in H_{j \mu}$ then $xy \in H_{i \mu}$.
\item[(ii)] Each idempotent is a left identity in its $\gr$-class and dually a right identity in its $\gl$-class. In other words $e_{i \omega} s = s$ for all $s \in R_i$, and  $s e_{i \omega} = s$ for all $s \in L_\omega$.
\item[(iii)] For all $i,j \in I$ and $\omega, \mu \in \Omega$ we have $H_{i \omega} \cong H_{j \mu}$.
\end{enumerate}
\end{proposition}

The Rees theorem (see \cite[Section~3.2]{Howie1995} or originally \cite{Rees1940}) characterises completely simple semigroups as Rees matrix semigroups over groups. We make use of Rees's result in Section~\ref{sec_FP1}.  

Let us fix some notation that will remain in force throughout this section.
Let $U$ be a completely simple semigroup
with set of $\gr$-, $\gl$- and $\gh$-classes
$
\{ R_i : i \in I \}$,
$\{ L_\omega : \omega \in \Omega \}$ and
$\{ H_{i \omega} : i \in I, \omega \in \Omega \}
$, respectively.
We assume throughout this section that $U$ has only finitely many left ideals, which is equivalent to saying that the set $\Omega$ is finite. 
We suppose that the index sets $I$ and $\Omega$ each contain the distinguished symbol $1$ and let
$R= R_1$, $L = L_1$ and $H = H_{1,1} = R_1 \cap L_1 = R \cap L$.
Define $S = U^1$ (the completely simple semigroup $U$ with an identity adjoined) and set $T = L \cup \{ 1 \} \subseteq S$ which is a submonoid of $S$.
Note that $L$ is a completely simple semigroup with underlying group $H$. 
In fact, $L$ is a left group. 
We have the inclusions
$H \leq T \leq S$. This allows us to break down the problem of relating $S$ with $H$ into two stages: first we relate $S$  and $T$ (in Proposition~\ref{prop_SiffT}), and then we relate $T$ and $H$ (in Proposition~\ref{prop_TiffH}). 
Theorem~\ref{thm_completelysimplemain} then follows by combining Propositions~\ref{prop_SiffT} and \ref{prop_TiffH}. 

It is important to observe that with the above definitions $\ZT$ is \emph{not} free when viewed as a left $\ZH$-module, and also $\ZS$ is \emph{not} free when viewed as a left $\ZT$-module.

The rest of this section will be devoted to the problem of relating free left resolutions for $S$ with those for $T$, with the main result being Proposition~\ref{prop_hash}, from which Proposition~\ref{prop_SiffT} is a consequence. 

Let $e \in H$ be the idempotent in the $\gh$-class $H$, and set let
$F$ denote the set of all other idempotents in $R$ i.e. $F = \{ E(U) \setminus \{ e \} \} \cap R$.
Note that $F$ is finite since the index set $\Omega$ is assumed to be finite.
For $y \in U$ we use $L_y$ to denote the $\gl$-class of $U$ containing $y$. Similarly we define $R_y$ and $H_y$. 

The following self-evident fact will be used repeatedly throughout the section. 

\begin{lemma}
\label{lem_decomposition}
Every element $\lambda \in \ZS$ may be written uniquely in the form
$$
\lambda = \lambda^{(1)} + \sum_{f \in F} \lambda^{(f)}
$$
where $\lambda^{(1)} \in \ZT$ and $\lambda^{(f)} \in \ZL_f$ for all $f \in F$.
\end{lemma}

The general observation given in the next lemma will be important to us. 

\begin{lemma}
\label{lem_changeoftype}
Let $A$ be a left $\ZS$-module. Then $A$ is finitely generated as a left $\ZS$-module if and only if $A$ is finitely generated as a left $\ZT$-module.
\end{lemma}
\begin{proof}
For the non-trivial direction of the proof let $X$ 
be a finite generating set for $A$ as a left $\ZS$-module. 
For every $\lambda \in \ZS$, decomposing as in Lemma~\ref{lem_decomposition}, we obtain
\[
\lambda
= \lambda^{(1)} + \sum_{f \in F} \lambda^{(f)}
= \lambda^{(1)} + \sum_{f \in F} \lambda^{(f)} ef
\]
where $\lambda^{(1)} \in \ZT$, $\lambda^{(f)} \in \ZL_f$ $(f \in F)$ and by Proposition~\ref{prop_compsimpleprops} for all $f \in F$, 
$\lambda^{(f)} = \lambda^{(f)} f =
\lambda^{(f)} e f$ and $\lambda^{(f)} e \in \ZL \subseteq \ZT$.
It follows that for all $x \in X$ and $\lambda \in \ZS$:
\[
\lambda x =
\lambda^{(1)} x + \sum_{f \in F} (\lambda^{(f)} e) (f x)
\in \lb X \cup FX \rb_{\ZT}.
\]
Therefore
\[
X \cup F X = X \cup \{ f x : f \in F, x \in X \} \subseteq A
\]
is a generating set for $A$ regarded as a left $\ZT$-module, where $X \cup FX$ is finite since $X$ and $F$ are both finite.  
\end{proof}

For every partial free left resolution of $S$ we shall associate a partial free left resolution of $T$ and mappings $\theta$ and $\phi$ relating the two resolutions. Let
\[
\mathcal{A}:
A_n \xrightarrow{\partial_n}
A_{n-1} \xrightarrow{\partial_{n-1}}
\cdots \xrightarrow{\partial_2} A_1
\xrightarrow{\partial_1} A_0 \xrightarrow{\partial_0} \mathbb{Z} \rightarrow 0
\]
be a partial free resolution of the trivial left $\ZS$-module $\mathbb{Z}$, where 
$A_0 = \ZS$, $\partial_0 = \epsilon_S$ is  the standard augmentation, $\ker \epsilon_S = \lb X_0 \rb_\ZT$, and for $j = 1, \ldots, n$
\[
A_j = \bigoplus_{x \in X_{j-1}} \ZS [x]
\]
where $\ker(\partial_{j-1}) = \lb X_{j-1} \rb_{\ZT}$ and $\partial_j: A_j \rightarrow A_{j-1}$ is the left $\ZS$-module homomorphism extending $[x] \mapsto x$ $(x \in X_{j-1})$.
Note here that each $X_{j-1}$ has been chosen so that $\lb X_{j-1} \rb_\ZT = \ker(\partial_{j-1})$ and not just $\lb X_{j-1} \rb_\ZS = \ker(\partial_{j-1})$. By Lemma~\ref{lem_Schanuel} and Lemma~\ref{lem_changeoftype} it follows that if $S$ is of type left-$\FPn$ then such a partial resolution $\mathcal{A}$ exists with $|X_i| < \infty$ for all $0 \leq i \leq n-1$.   

Using $\mathcal{A}$, our aim is to construct a resolution $\mathcal{B}$ for $T$. Define 
\[
B_0 = \ZT[e] \oplus \bigoplus_{f \in F} \ZT[f],
\]
and, using the natural decomposition given in Lemma~\ref{lem_decomposition}, define a mapping $\theta: A_0 \rightarrow B_0$ by 
\[
\theta(\lambda) = \lambda^{(1)} [e] + \sum_{f \in F} \lambda^{(f)} e [f] \quad 
(\lambda \in A_0 = \ZS).
\]
Then for $1 \leq m \leq n$ define
\[
B_m =
( \bigoplus_{x \in X_{m-1}} \ZT [x] )
\oplus
( \bigoplus_{i = m-1, \ldots, 0 \atop x \in X_{i}, \; f \in F } \ZT [f,x] )
\oplus
( \bigoplus_{f \in F } \ZT [f] ),
\]
and a mapping $\theta: A_m \rightarrow B_m$ given by
\[
\theta( \sum_{x \in X_{m-1}} \lambda_x [x])
=
\sum_{x \in X_{m-1}} \lambda_x^{(1)} [x] +
\sum_{x \in X_{m-1} \atop f \in F}
\lambda_x^{(f)} e [f,x]
\]
where $\lambda_x \in \ZS$ for $x \in X_{m-1}$.
The fact that $\theta(A_m)$ is a subset of $B_m$ follows from the definition of $\theta$ along with Proposition~\ref{prop_compsimpleprops}.
We view $\theta$ as a mapping with domain $\dotcup_{0 \leq i \leq n} A_i$ such that for every $1 \leq m \leq n$ the restriction of $\theta$ to $A_m$ maps $A_m$ to $B_m$. 

We may now state the main result of this section which relates free resolutions of $S$ with free resolutions of $T$.

\begin{proposition}
\label{prop_hash}
Let $U$ be a completely simple semigroup with finitely many left ideals, let $L$ be an $\gl$-class of $U$, and set $S = U^1$ and $T = L^1$. 
Let
\[
\mathcal{A}:
A_n \xrightarrow{\partial_n}
A_{n-1} \xrightarrow{\partial_{n-1}}
\cdots \xrightarrow{\partial_2} 
A_1
\xrightarrow{\partial_1} 
A_0 \xrightarrow{\partial_0 = \epsilon_S} \mathbb{Z} \rightarrow 0
\]
be a partial free resolution of the trivial left $\ZS$-module $\mathbb{Z}$ where $A_0 = \ZS$, $\partial_0 = \epsilon_S$ is the standard augmentation, $\ker \epsilon_S = \lb X_0 \rb_\ZT$ and for $1 \leq j \leq n$ 
\[
A_j = \bigoplus_{x \in X_{j-1}} \ZS [x]
\]
where $\ker(\partial_{j-1}) = \lb X_{j-1} \rb_{\ZT}$ and $\partial_j: A_j \rightarrow A_{j-1}$ is the left $\ZS$-module homomorphism extending $[x] \mapsto x$ $(x \in X_{j-1})$.
Then with the above notation: 
\[
\mathcal{B}:
B_n \xrightarrow{\partial_n'}
B_{n-1} \xrightarrow{\partial_{n-1}'}
\cdots \xrightarrow{\partial_2'}
B_1 \xrightarrow{\partial_1'}
B_0 \xrightarrow{\partial_0'} \mathbb{Z} \rightarrow 0
\]
is a partial free resolution of the trivial left $\ZT$-module $\mathbb{Z}$ where
$\partial_0' : B_0 \rightarrow \mathbb{Z}$ is the left $\ZS$-module homomorphism extending:
\[
\partial_0'(\lambda_e [e] + \sum_{f \in F} \lambda_f [f]) = \epsilon_T(\lambda_e) + \sum_{f \in F} \epsilon_T(\lambda_f),
\]
where $\epsilon_T$ is the standard augmentation, 
and for $1 \leq m \leq n$,
$\partial_m': B_m \rightarrow B_{m-1}$ is the left $\ZT$-module homomorphism extending: $\partial_m' ([x]) = \theta(x)$ $(x \in X_{m-1})$, and for $f \in F$ and $x \in X_i$,
\[
\partial_m'([f,x]) = \begin{cases}
\theta(fx) & \mbox{if $i = m-1$} \\
(1-e)[f,x] & \mbox{if $i=m-2, m-4, \ldots$} \\
e[f,x] & \mbox{if $i=m-3, m-5, \ldots$}
\end{cases}
\]
and
\[
\partial_m' ([f]) = \begin{cases}
e[f] & \mbox{if $m$ is even} \\
(1-e)[f] & \mbox{if $m$ is odd}.
\end{cases}
\]
Furthermore $\ker \partial_0$ is finitely generated if and only if $\ker \partial_0'$ is finitely generated; and
if $A_0, A_1, \ldots, A_n$ are all finitely generated then $B_0, B_1, \ldots, B_n$ are finitely generated as well, in which case $\ker \partial_n$ is finitely generated if and only if $\ker \partial_n'$ is finitely generated.
\end{proposition}

\begin{remark}
\label{remark_theta}
Note that the definition of $\partial_m'$ given in Proposition~\ref{prop_hash} makes sense since for $x \in X_{m-1}$ and $f \in F$ we have $x, fx \in A_{m-1}$ which by definition of $\theta$ implies 
\begin{equation}
\label{eqn_remark3}
\theta(x), \; \theta(fx) \in
\bigoplus_{x \in X_{m-2}} \ZT[x]
\oplus
\bigoplus_{x \in X_{m-2} \atop f \in F} \ZT[f,x]
\subseteq B_{m-1}.
\end{equation}
\end{remark}

Before embarking on the proof of Proposition~\ref{prop_hash} we shall give one of its important consequences. 

\begin{proposition}
\label{prop_SiffT}
Let $U$ be a completely simple semigroup and let $L$ be a left ideal of $U$. If $U$ has finitely many left ideals then $U$ is of type left-$\FPn$ if and only if $L$ is of type left-$\FPn$. 
\end{proposition}
\begin{proof}
This is equivalent (in the above notation) to proving that if $F$ is finite then $S = U^1$ is of type left-$\FP_n$ if and only if $T = L^1$ is of type left-$\FP_n$.
We prove the result by induction on $n$. The base case $n=1$ is easy to verify directly. Alternatively, it follows from a more general result that we prove below (see Theorem~\ref{thm_fp1compsimple} in Section~\ref{sec_FP1}). Now let $n>1$ and assume inductively that the result holds for values strictly less than $n$. 

If $S$ is of type left-$\FPn$ then by Lemma~\ref{lem_Schanuel} and Lemma~\ref{lem_changeoftype} the partial resolution $\mathcal{A}$ given in the statement of Proposition~\ref{prop_hash} may be chosen so that $X_i$ is finite for all $0 \leq i \leq n-1$. By Proposition~\ref{prop_hash}, $\mathcal{B}$ is then a partial free resolution of the trivial left $\ZT$-module $\mathbb{Z}$ with each $B_i$ free of finite rank. Therefore $T$ is of type left-$\FPn$. 

Conversely suppose that $T$ is of type left-${\rm FP}\sb n$. In particular $T$ is of type left-$\FP_{n-1}$ which by induction implies that $S$ is of type left-$\FP_{n-1}$. By Lemma~\ref{lem_Schanuel} and Proposition~\ref{prop_hash}, using the definitions given in Proposition~\ref{prop_hash}, we have that 
\[
\begin{array}{ccccccccccc}
\mathcal{A}: A_{n-1}  & \xrightarrow{\partial_{n-1}} &
A_{n-2}  & \xrightarrow{\partial_{n-2}} &
\cdots   & \xrightarrow{\partial_{1}} &
A_0      & \xrightarrow{\partial_0 = \epsilon_S} &
\mathbb{Z} & \rightarrow & 0, \\
& &
& & & &
& &
& & \\
\mathcal{B}: B_{n-1}  & \xrightarrow{\partial_{n-1}'} &
B_{n-2}  & \xrightarrow{\partial_{n-2}'} &
\cdots   & \xrightarrow{\partial_{1}'} &
B_0      & \xrightarrow{\phantom{=} \partial_0' \phantom{\epsilon_T}} &
\mathbb{Z} & \rightarrow & 0
\end{array}
\]
are both partial free resolutions with $A_j$ and $B_j$ finitely generated 
 for $j=0, 1, \ldots, n-1$. By Lemma~\ref{lem_Schanuel} since $T$ is of type left-${\rm FP}\sb n$ it follows that $\ker \partial_{n-1}'$ is finitely generated which by the last clause in the statement of Proposition~\ref{prop_hash} implies that $\ker \partial_{n-1}$ is finitely generated (both as a left $\ZS$-module and a left $\ZT$-module, by Lemma~\ref{lem_changeoftype}). 
 It now follows from Lemma~\ref{lem_Schanuel} that $S$ is of type left-${\rm FP}\sb n$.
\end{proof}

We now work through several technical lemmas which will then be utilised at the end of the section where we prove Proposition~\ref{prop_hash}.

First we define a mapping $\phi: B_m \rightarrow A_m$ which taken together with $\theta: A_m \rightarrow B_m$ (defined above) will help clarify the relationship between $\mathcal{A}$ and $\mathcal{B}$.

Define $\phi: B_0 \rightarrow A_0$ by:
\[
\phi(\lambda_e [e] + \sum_{f \in F} \lambda_f [f]) = \lambda_e + \sum_{f \in F} \lambda_f f,
\]
and for $1 \leq m \leq n$ define $\phi: B_m \rightarrow A_m$ by:
\begin{eqnarray*}
& & \phi(
\sum_{x \in X_{m-1}} \lambda_x [x] +
\sum_{i = m-1, \ldots, 0 \atop x \in X_{i}, f \in F } \lambda_{f,x} [f,x] +
\sum_{f \in F} \lambda_f [f]
) \\
& = &
\sum_{x \in X_{m-1}} \lambda_x [x] +
\sum_{x \in X_{m-1} \atop f \in F} \lambda_{f,x} f[x].
\end{eqnarray*}

In a similar way as for $\theta$ above, we view $\phi$ as a map with domain $\dotcup_{0 \leq i \leq n} B_i$ where the restriction of $\phi$ to $B_m$ maps $B_m$ to $A_m$. 
The relationship between the mappings thus far defined is illustrated below: 
\[
\begin{array}{ccccccccccc}
A_n      & \xrightarrow{\partial_n} &
A_{n-1}  & \xrightarrow{\partial_{n-1}} &
\cdots & \xrightarrow{\partial_{1}} &
A_0      & \xrightarrow{\partial_0 = \epsilon_S} &
\mathbb{Z} & \rightarrow & 0 \\
\theta \downarrow \uparrow \phi & &
\theta \downarrow \uparrow \phi & & 
& &
\theta \downarrow \uparrow \phi & &
& & \\
B_n      & \xrightarrow{\partial_n'} &
B_{n-1}  & \xrightarrow{\partial_{n-1}'} &
\cdots & \xrightarrow{\partial_{1}'} &
B_0      & \xrightarrow{\phantom{= } \partial_0' \phantom{\epsilon_T}} &
\mathbb{Z} & \rightarrow & 0.
\end{array}
\]

The next lemma tells us that $\theta$ behaves well with respect to addition and the action of $\ZT$.

\begin{lemma}
\label{lem_propoftheta}
For all $0\leq m \leq n$ the mapping $\theta: A_m \rightarrow B_m$ is a homomorphism of abelian groups and commutes with the action of $\ZT$.
That is, for all $\lambda \in \ZT$ and $\alpha \in A_m$ we have 
$
\lambda \cdot \theta(\alpha) = \theta( \lambda \cdot \alpha).
$
\end{lemma}
\begin{proof}
For all $\lambda, \mu \in \ZS$ and $f \in F$, decomposing as in Lemma~\ref{lem_decomposition}, it is easily seen that
\[
(\lambda + \mu)^{(1)} = \lambda^{(1)} + \mu^{(1)} \quad \mbox{and} \quad
(\lambda + \mu)^{(f)} = \lambda^{(f)} + \mu^{(f)}.
\]
From this and the definition of $\theta$ it follows that
\[
\theta(\alpha + \beta) = \theta(\alpha) + \theta(\beta)
\]
for all $\alpha, \beta \in A_m$. Hence $\theta$ is a homomorphism of abelian groups. 

By Proposition~\ref{prop_compsimpleprops}, under its action on $S$ by
left multiplication, $T$ stabilises itself and each of the sets $L_f \ (f \in F)$ setwise.
It follows that for all $\mu \in \ZS$, $\lambda \in \ZT$ and $f \in F$:
\[
(\lambda \cdot \mu)^{(1)} =  \lambda \cdot \mu^{(1)} 
\quad \mbox{and} \quad
(\lambda \cdot \mu)^{(f)} = \lambda \cdot \mu^{(f)}.  
\]
From this and the definition of $\theta$ we conclude that
$
\lambda \cdot \theta(\alpha) = \theta( \lambda \cdot \alpha)
$
for all $\lambda \in \ZT$ and $\alpha \in A_m$.
\end{proof}

The map $\phi$ is equally well behaved. 

\begin{lemma}\label{lem_propofphi}
For all $0 \leq m \leq n$
the mapping $\phi: B_m \rightarrow A_m$ is a homomorphism of abelian groups and commutes with the action of $\ZT$.
That is, for all $\mu \in \ZT$ and $\alpha \in B_m$ we have
$
\mu \cdot \phi(\alpha) = \phi( \mu \cdot \alpha).
$
\end{lemma}
\begin{proof}
This follows easily from the definition of $\phi$. 
\end{proof}

\begin{lemma}
\label{lem_identitymap}
For $0 \leq m \leq n$
the composition $\phi \theta: A_m \rightarrow A_m$ is the identity map on $A_m$. In particular, $\theta$ is injective.
\end{lemma}
\begin{proof}
Suppose that $m \geq 1$, the case $m=0$ may be handled similarly. Let 
$\sum_{x \in X_{m-1}} \lambda_x [x] \in A_m$, where $\lambda_x \in \ZS$ for $x \in X_{m-1}$. Then from the definitions of $\phi$ and $\theta$ we have
\begin{eqnarray*}
\phi \theta ( \sum_{x \in X_{m-1}} \lambda_x [x] )
& = &
\phi(
\sum_{x \in X_{m-1}} \lambda_x^{(1)} [x] +
\sum_{x \in X_{m-1} \atop f \in F}
\lambda_x^{(f)} e [f,x]) \\
& = &
\sum_{x \in X_{m-1}} \lambda_x^{(1)} [x] +
\sum_{x \in X_{m-1} \atop f \in F}
\lambda_x^{(f)} e f [x]
= \sum_{x \in X_{m-1}} \lambda_x [x]
\end{eqnarray*}
since, by Proposition~\ref{prop_compsimpleprops}, for $f \in F$ and $x \in X_{m-1}$ we have
$
\lambda_x^{(f)} ef= \lambda_x^{(f)} f = \lambda_x^{(f)},
$
because $\lambda_x^{(f)} \in L_f$ and $f$ is a right identity in its $\gl$-class.
\end{proof}

\begin{lemma}
\label{phikerneltokernel}
We have
$\phi(\ker \partial_m') \subseteq \ker \partial_m$
for all $0 \leq m \leq n$.
\end{lemma}
\begin{proof}
The fact that $\phi(\ker \partial_0')$ is a subset of $\ker \partial_0$ is an easy consequence of the definition
of $\phi: B_0 \rightarrow A_0$. This 
deals with the case $m=0$.
Now let $l \in \ker \partial_m'$ for some $1 \leq m \leq n$. Say:
\[
l =
\sum_{x \in X_{m-1}} \lambda_x [x] +
\sum_{i = m-1, \ldots, 0, \atop x \in X_{i}, f \in F } \lambda_{f,x} [f,x] +
\sum_{f \in F} \lambda_f [f],
\]
where each of the terms $\lambda_x, \lambda_{f,x}$ and $\lambda_f$ belongs to $\ZT$. 
Suppose that $m$ is odd, the case when $m$ is even is dealt with similarly. Then
since $l \in \ker \partial_m'$ applying $\partial_m'$ gives 
\begin{align}
\label{eqn_triangle2}
\sum_{x \in X_{m-1}} \lambda_x \theta(x) +
\sum_{x \in X_{m-1} \atop f \in F} \lambda_{f,x} \theta(fx) +  
\sum_{i = m-2, m-4, \ldots  \atop x \in X_{i}, f \in F } \lambda_{f,x} (1-e) [f,x]   
\quad \quad \quad \quad \notag
\\
+ \sum_{i = m-3, m-5, \ldots  \atop x \in X_{i}, f \in F } \lambda_{f,x} e [f,x] +
\sum_{f \in F} \lambda_f (1-e) [f] = \partial_m' (l) = 0. 
\end{align}
From this along with \eqref{eqn_remark3} and the definition of $\theta: A_{m-1} \rightarrow B_{m-1}$ we deduce:
\begin{eqnarray*}
\lambda_{f,x} e = 0
& &
\mbox{for $f \in F, \; x \in X_i \ (i = m-3, m-5, \ldots )$} \\
\lambda_{f,x} (1-e) = 0
& &
\mbox{for $f \in F, \; x \in X_i \ (i = m-4, m-6, \ldots )$} \\
\lambda_f (1-e) = 0
& &
\mbox{for $f \in F$},
\end{eqnarray*}
and hence equation \eqref{eqn_triangle2} becomes:
\begin{equation}
\label{eqn_star2}
\sum_{x \in X_{m-1}} \lambda_x \theta(x) +
\sum_{x \in X_{m-1} \atop f \in F} \lambda_{f,x} \theta(fx)
+
\sum_{x \in X_{m-2} \atop f \in F} \lambda_{f,x} (1-e) [f,x]  =0.
\end{equation}
(When $m=1$ the third term here is $\sum_{f \in F} \lambda_f (1-e) [f]$ but the rest of the argument follows the same lines as below.)
It follows from equation \eqref{eqn_star2} and Lemma~\ref{lem_propoftheta} that 
\[
\sum_{x \in X_{m-2} \atop f \in F} \lambda_{f,x} (1-e) [f,x] \in \theta(A_{m-1}) \subseteq B_{m-1},
\]
which from the definition of $\theta$ implies that $\lambda_{f,x}(1-e) \in \ZL$, for all $f \in F$ and $x \in X_{m-2}$. Therefore $\lambda_{f,x}(1-e) = \lambda_{f,x}(1-e)e = 0$ since $e$ is a right identity in the $\gl$-class $L$. 
Now substituting this back into equation \eqref{eqn_star2} and applying Lemma~\ref{lem_propoftheta} gives
\begin{eqnarray*}
\theta(\partial_m (\phi(l)) )
&  = &
\theta (\sum_{x \in X_{m-1}} \lambda_x x +
\sum_{x \in X_{m-1} \atop f \in F} \lambda_{f,x} fx )
\\
& = &
\sum_{x \in X_{m-1}} \lambda_x \theta(x) +
\sum_{x \in X_{m-1} \atop f \in F} \lambda_{f,x} \theta(fx) = 0.
\end{eqnarray*}
But by Lemma~\ref{lem_identitymap}, $\theta$ is injective and so since $\theta(0) = 0 = \theta(\partial_m(\phi(l)))$ we conclude $\partial_m (\phi(l)) = 0$.
\end{proof}

\begin{lemma}
\label{lem_thetakerneltokernel}
We have
$\theta(\ker \partial_m) \subseteq \ker \partial_m'$
for all $0 \leq m \leq n$.
\end{lemma}
\begin{proof}
It is an immediate consequence of the definition of $\theta: A_0 \rightarrow B_0$ that $\theta( \ker \partial_0 ) \subseteq \ker \partial_0'$. This shows that the result holds when $m=0$. 
Next let $1 \leq m \leq n$ and take an arbitrary element $\alpha = \sum_{x \in X_{m-1} } \lambda_x [x] \in \ker \partial_m$. Then from the definitions and Lemma~\ref{lem_propoftheta} we have:
\begin{eqnarray*}
\partial_m' (\theta( \sum_{x \in X_{m-1} } \lambda_x [x] ) )
& = &
\partial_m' (
\sum_{x \in X_{m-1}} \lambda_x^{(1)} [x] +
\sum_{x \in X_{m-1} \atop f \in F}
\lambda_x^{(f)} e [f,x]
)
\\
& =&
\sum_{x \in X_{m-1}} \lambda_x^{(1)} \theta(x) +
\sum_{x \in X_{m-1} \atop f \in F}
\lambda_x^{(f)} e \theta(fx)
\\
& = &
\theta(
\sum_{x \in X_{m-1}} \lambda_x^{(1)} x +
\sum_{x \in X_{m-1} \atop f \in F}
\lambda_x^{(f)} e fx
)
\\
& = &
\theta(
\sum_{x \in X_{m-1}} \lambda_x x
)
 = 
\theta(\partial_m (\alpha)) = \theta(0) = 0,
\end{eqnarray*}
where $\lambda_x^{(f)} ef = \lambda_x^{(f)} f = \lambda_x^{(f)}$ by Proposition~\ref{prop_compsimpleprops} since $\lambda_x^{(f)} \in \ZL_f$.
\end{proof}

The next result relates generating sets of kernels from the sequence $\mathcal{B}$ with generating sets of kernels from the sequence $\mathcal{A}$. 

\begin{lemma}\label{lem_StoT2}
Let $Y$ be a subset of $\ker \partial_m'$ for some $0 \leq m \leq n$. If $\lb Y \rb_\ZT = \ker \partial_m'$ then $\lb \phi(Y) \rb_\ZT = \ker \partial_m$.
\end{lemma}
\begin{proof}
It follows from Lemma~\ref{phikerneltokernel} that $\phi(Y) \subseteq \ker \partial_m$.  
Let $\alpha \in \ker \partial_m $ be arbitrary. Since $\phi \theta$ is the identity on $A_m$ by Lemma~\ref{lem_identitymap}, we have $\alpha = \phi \theta (\alpha)$ where $\theta (\alpha) \in \ker \partial_m' = \lb Y \rb_\ZT$ by Lemma~\ref{lem_thetakerneltokernel}.
Write:
\[
\theta (\alpha) = \sum_{y \in Y} \gamma_y y, \quad (\gamma_y \in \ZT).
\]
Then by Lemma~\ref{lem_propofphi}:
\[
\alpha = \phi \theta (\alpha)
 =
\phi(
\sum_{y \in Y} \gamma_y y
)
=
\sum_{y \in Y}
\gamma_y \phi(y) \in
\lb \phi(Y) \rb_\ZT.
\]
Since $\alpha$ was arbitrary it follows that $\lb \phi(Y) \rb_\ZT = \ker \partial_m$. 
\end{proof}

We are now in a position to prove the main result of this section. 

\begin{proof}[Proof of Proposition~\ref{prop_hash}]
For $0 \leq m \leq n$, we define the following subset $Y_m$ of $B_m$:
\begin{eqnarray*}
Y_m = \theta(X_m)
& \cup &
\{
(1-e)[f,x] : f \in F, \ x \in X_i \ (i = m-1, m-3, \ldots )
\} \\
& \cup &
\{
e[f,x] : f \in F, \ x \in X_i \ (i = m-2, m-4, \ldots )
\} \\
& \cup &
Q,
\end{eqnarray*}
where
\[
Q =
\begin{cases}
\{ (1-e)[f] : f \in F \} & \mbox{if $m$ is even} \\
\{ e[f] : f \in F \} & \mbox{if $m$ is odd}.
\end{cases}
\]

\bigskip

\noindent \textbf{Claim.} For $0 \leq m \leq n$ we have $Y_m \subseteq \ker \partial_m'$ and $\lb Y_m \rb_\ZT = \ker \partial_m'$.

\smallskip

\begin{proof}[Proof of Claim.]
Suppose that $m$ is odd, the case when $m$ is even is similar. In particular the case $m=0$ may be dealt with using a similar argument.

Since $X_m \subseteq \ker \partial_m$ it follows from Lemma~\ref{lem_thetakerneltokernel} that $\theta(X_m) \subseteq \ker \partial_m'$.
It then rapidly follows from the definition of $\partial_m'$ and Lemma~\ref{lem_propoftheta} that $Y_m \subseteq \ker \partial_m'$.

To see that $\lb Y_m \rb_\ZT = \ker \partial_m'$ let $l \in \ker \partial_m'$ be arbitrary, say
\[
l =
\sum_{x \in X_{m-1}} \lambda_x [x] +
\sum_{i = m-1, \ldots, 0, \atop x \in X_{i}, f \in F } \lambda_{f,x} [f,x] +
\sum_{f \in F} \lambda_f [f],
\]
where each of the terms $\lambda_x, \lambda_{f,x}$ and $\lambda_f$ belongs to $\ZT$. 
Exactly as in the proof of Lemma~\ref{phikerneltokernel} from $\partial_m'(l)=0$ we deduce:
\begin{eqnarray*}
\lambda_{f,x} e = 0
& &
\mbox{for $f \in F \setminus \{e \}, \; x \in X_i \ (i = m-3, \; m-5, \ldots )$} \\
\lambda_{f,x} (1-e) = 0
& &
\mbox{for $f \in F \setminus \{e \}, \; x \in X_i \ (i = m-2, \; m-4, \ldots )$} \\
\lambda_f (1-e) = 0.
\end{eqnarray*}
It then follows from the definition of $Y_m$, and since $m$ is odd, that
\begin{eqnarray*}
& &
\sum_{i = m-2, m-3 \ldots, 0 \atop
x \in X_{i}, f \in F} \lambda_{f,x} [f,x]
+
\sum_{f \in F} \lambda_f [f] \\
& = &
\sum_{i = m-2, m-4, \ldots \atop
x \in X_{i}, f \in F} \lambda_{f,x} e [f,x]
+
\sum_{i = m-3, m-5, \ldots \atop
x \in X_{i}, f \in F} \lambda_{f,x} (1-e) [f,x]
+
\sum_{f \in F} \lambda_f e [f] \\
& \in & \lb Y_m \rb_\ZT.
\end{eqnarray*}
In other words
$l-l_1 \in \lb Y_m \rb_\ZT$ where
\[
l_1 =
\sum_{x \in X_{m-1}} \lambda_x [x] +
\sum_{x \in X_{m-1} \atop f \in F } \lambda_{f,x} [f,x] \in \ker \partial_m'.
\]
Since $l_1 \in \ker \partial_m'$, by Lemma~\ref{phikerneltokernel} we have $\phi(l_1) \in \ker \partial_m = \lb X_m \rb_\ZT$, and applying
Lemma~\ref{lem_propoftheta}:
\[
\theta(\phi(l_1)) 
=
\sum_{x \in X_{m-1}} \lambda_x [x] +
\sum_{x \in X_{m-1} \atop f \in F } \lambda_{f,x} e [f,x]
\in \lb \theta(X_m) \rb_\ZT \subseteq \lb Y_m \rb_\ZT.
\]
We conclude, again by inspection of $Y_m$, that 
\[
l_1
=
\theta(\phi(l_1))
+
\sum_{x \in X_{m-1} \atop f \in F } \lambda_{f,x} (1-e) [f,x]
\in \lb Y_m \rb_\ZT,
\]
and hence $l = (l-l_1) + l_1 \in \lb Y_m \rb_\ZT$, completing the proof of the claim. 
\end{proof}
Returning to the proof of Proposition~\ref{prop_hash}, 
by definition each of $B_0$, $B_1$, $\ldots, B_n$ is a free left $\ZT$-module. Next we must show that $\mathcal{B}$ is exact. It is an immediate consequence of the definitions that $\partial_0'(B_0) = \mathbb{Z}$. Now let $1 \leq j \leq n-1$ and consider $\im \partial_{j+1}'$.
For all $x \in X_j \subseteq \ker(\partial_j)$ and $f \in F$ we have $fx \in \ker(\partial_j)$ and so by Lemma~\ref{lem_thetakerneltokernel} and the claim:
\[
\theta(fx) \in \theta(\ker(\partial_j)) \subseteq \ker \partial_j' = \lb Y_j \rb_\ZT.
\]
Using this observation, it is then easy to verify that $\im \partial_{j+1}' \subseteq \lb Y_j \rb_\ZT$. 
From the definition of $\partial_{j+1}'$ we see that $Y_j \subseteq \im \partial_{j+1}'$. Therefore, applying the above claim we conclude that $\mathcal{B}$ is exact. 

The last clauses in the statement of Proposition~\ref{prop_hash} follow from the claim and Lemma~\ref{lem_StoT2}. 
\end{proof}

\section{Resolutions for Completely Simple Semigroups II}
\label{sec_resolutions2}

In this section we finish off the proof of Theorem~\ref{thm_completelysimplemain}.
In Section~\ref{sec_resolutions} above, 
in the notation of that section, 
we saw how to pass between free resolutions of $S$ and free resolutions of $T$.
We now go on to consider the relationship between resolutions of $T$ and those of $H$
with the aim of proving the following result. 

\begin{proposition}
\label{prop_TiffH}
Let $L$ be a left group and let $H$ be a maximal subgroup of $L$. Then $L$ is of type left-$\FPn$ if and only if $L$ has finitely many idempotents, and $H$ is of type $\FPn$.  
\end{proposition}

Recall that a left group is the direct product of a left zero semigroup and a group. 
We note that a recent result of Kobayashi \cite[Corollary~2.7]{Kobayashi2007} is exactly Proposition~\ref{prop_TiffH} in the special case $n=1$.  

Note that Proposition~\ref{prop_TiffH} does not simply follow from the dual of Proposition~\ref{prop_SiffT}, since all statements are about \emph{left} $\ZS$-modules. 

Proposition~\ref{prop_TiffH} will follow from Proposition~\ref{prop_diamond} and Theorem~\ref{thm_fp1compsimple} which will be proved below. 
We continue using the same notation that was introduced above in Section~\ref{sec_resolutions}, with the exception of $F$ which will be used to denote a different set of idempotents from before.   
So $T = L \cup \{ 1 \} \leq S$ where $L = L_1$ is a completely simple semigroup with a single $\gl$-class.  Fix an idempotent $e$ in $L$, 
let $F = E(L)$ and set $H = H_e$.
For $f \in F$ we use $H_f$ to denote the $\gh$-class of $f$. 
Note that now $F$ is the set of idempotents of an $\gl$-class, while in the previous section above $F$ was used to denote a set of $\gr$-related idempotents. 

One direction of Proposition~\ref{prop_TiffH} is straightforward.
Suppose that $T$ is of type left-$\FP_n$. Then $H$ is of type left-$\FPn$ by Theorem~\ref{thm_MinimalMaximalSubgroups}. 
Also, since $T$ is of type left-$\FP_1$, it follows from Theorem~\ref{thm_fp1compsimple} below (see also
\cite[Theorem~2.6]{Kobayashi2007}) that $L$ has finitely many idempotents.

The main result of this section is the following which, when combined with the previous paragraph, has  
Proposition~\ref{prop_TiffH} as a consequence.  

\begin{proposition}
\label{prop_diamond}
Let $L$ be a left group, $H$ be a maximal subgroup of $L$ with identity $e$, and set $F=E(L)$ and $T=L^1$. 
Let
\[
\mathcal{A}:
A_n \xrightarrow{\partial_n}
A_{n-1} \xrightarrow{\partial_{n-1}}
\cdots \xrightarrow{\partial_2} A_1
\xrightarrow{\partial_1} A_0 \xrightarrow{\partial_0 = \epsilon_H} \mathbb{Z} \rightarrow 0
\]
be a partial free resolution of the trivial left $\ZH$-module $\mathbb{Z}$ where $\partial_0 = \epsilon_H$ is the standard augmentation, $\ker \epsilon_H = \lb X_0 \rb_\ZH$ and for $1 \leq m \leq n$
\[
A_m = \bigoplus_{x \in X_{m-1}} \ZH [x]
\]
where $\ker(\partial_{m-1}) = \lb X_{m-1} \rb_{\ZH}$ and $\partial_m: A_m \rightarrow A_{m-1}$ is the left $\ZH$-module homomorphism extending $[x] \mapsto x$ $(x \in X_{m-1})$. Then
\[
\mathcal{B}:
B_n \xrightarrow{\partial_n'}
B_{n-1} \xrightarrow{\partial_{n-1}'}
\cdots \xrightarrow{\partial_2'}
B_1 \xrightarrow{\partial_1'}
B_0 \xrightarrow{\partial_0' = \epsilon_T} \mathbb{Z} \rightarrow 0
\]
is a partial free resolution of the trivial left $\ZT$-module $\mathbb{Z}$ where $B_0 = \ZT$, $\partial_0' = \epsilon_T$ is the standard augmentation, and for $1 \leq m \leq n$:
\[
B_m =
\bigoplus_{i = m-1, \ldots, 0 \atop x \in X_i}
\ZT[x]
\oplus
\bigoplus_{f \in F}
\ZT [f],
\]
where, viewing $A_m \subseteq B_m$ under the natural inclusion arising from $\ZH \subseteq \ZT$, $\partial_m' : B_m \rightarrow B_{m-1}$ is the left $\ZT$-module homomorphism extending:
\[
\partial_m'([x]) = \begin{cases}
x & \mbox{if $x \in X_{m-1}$} \\
(1-e)[x] & \mbox{if $x \in X_i$, $i = m-2, \; m-4, \ldots$} \\
e[x] & \mbox{if $x \in X_i$, $i = m-3, \; m-5, \ldots$} 
\end{cases}
\]
and
\[
\partial_m' ([f]) = \begin{cases}
(f-1) & \mbox{if $m=1$} \\
(f-1)[f] & \mbox{if $m$ odd and $m \neq 1$} \\
f[f] & \mbox{if $m$ even}.
\end{cases}
\]	
Moreover if $F$ is finite, and $A_0, A_1, \ldots, A_n$ are all finitely generated 
then
$B_0, B_1, \ldots, B_n$ are finitely generated. 
Thus if $H$ is of type left-$\FP_n$ then $T$ is of type left-$\FP_n$. 
\end{proposition}
\begin{proof}
Let $Y_0 = X_0 \cup \{ (1-f): f \in F \}$, and then for 
$1 \leq m \leq n$ define:
\begin{eqnarray*}
Y_m  = X_m
& \cup &
\{
(1-e)[x] : x \in X_{i} \ (i = m-1, m-3, \ldots )
\} \\
& \cup &
\{
e[x]: x \in X_{i} \ (i = m-2, m-4, \ldots )
\} \\
& \cup &
Q,
\end{eqnarray*}
where
\[
Q =
\begin{cases}
\{ (1-f)[f] : f \in F \} & \mbox{if $m$ is even} \\
\{ f[f] : f \in F \} & \mbox{if $m$ is odd}.
\end{cases}
\]

\bigskip

\noindent \textbf{Claim.} For $0 \leq m \leq n$ we have $Y_m \subseteq \ker \partial_m'$ and $\lb Y_m \rb_\ZT = \ker \partial_m'$.

\smallskip

\begin{proof}[Proof of Claim.]
We prove the result for $m$ odd. The proof in the case $m$ even (and in particular the case $m=0$) is similar.

First we must verify that $Y_m$ is a subset of $\ker \partial_m'$. 
From the definitions we see that $\partial_m$ is the restriction of $\partial_m'$ to $A_m \subseteq B_m$, and it follows from this that $X_m \subseteq \ker \partial_m'$ since $X_m \subseteq \ker \partial_m$. For $x \in X_{m-1}$ since $e$ is a left identity for $H$ we have
\[
\partial_m'((1-e)[x]) = (1-e)x = 0 \quad (x \in X_{m-1}). 
\]
It is then easily verified from the definition of $\partial_m'$ that the remaining members of $Y_m$ belong to $\ker \partial_m'$. This proves $Y_m \subseteq \ker \partial_m'$. 

To see that $\lb Y_m \rb_\ZT = \ker \partial_m'$, let $\alpha \in \ker \partial_m'$ be arbitrary, say:
\[
\alpha = 
\sum_{i = m-1, \ldots, 0 \atop x \in X_i}
\lambda_x [x]
+
\sum_{f \in F}
\lambda_f [f],
\]
where each $\lambda_x, \lambda_f \in \ZT$, and since $m$ is odd $\partial_m'(\alpha) = 0$ implies
\begin{eqnarray}
\label{eqn_TiffH1}
& & 
\sum_{x \in X_{m-1}} \lambda_x x
+
\sum_{i = m-2, m-4, \ldots \atop x \in X_i}
\lambda_x (1-e) [x] \notag
\quad \quad \quad 
\quad \quad \quad 
\quad \quad \quad 
\\
& & 
\quad \quad \quad 
\quad \quad \quad 
\quad \quad \quad 
+
\sum_{i = m-3, m-5, \ldots \atop x \in X_i}
\lambda_x e [x]
+
\sum_{f \in F}
\lambda_f (f-1) [f]
=
0.
\end{eqnarray}
Recalling that
\begin{equation}
\label{eqn_newstar}
X_{m-1} \subseteq \ker \partial_{m-1} \subseteq A_{m-1} = 
\bigoplus_{x \in X_{m-2}} \ZH [x],
\end{equation} 
it is then immediate from equation \eqref{eqn_TiffH1} that  
\begin{eqnarray*}
\lambda_{x} e = 0
& &
\mbox{for $x \in X_i \ (i = m-3, \; m-5, \ldots )$} \\
\lambda_{x} (1-e) = 0
& &
\mbox{for $x \in X_i \ (i = m-4, \; m-6, \ldots )$}
\\
\lambda_f (f-1) = 0 
& & 
\mbox{for $f \in F.$} 
\end{eqnarray*}
(When $m=1$ things work slightly differently, as we shall explain below.)
From these observations, along with the definition of $Y_m$, we deduce:
\begin{eqnarray*}
& &
\sum_{i = m-3, m-4, \ldots, 0 \atop x \in X_i}
\lambda_x [x]
+
\sum_{f \in F}
\lambda_f [f]
\\
& = &
\sum_{i = m-3, m-5, \ldots \atop x \in X_i}
\lambda_x (1-e) [x]
+
\sum_{i = m-4, m-6, \ldots \atop x \in X_i}
\lambda_x e [x]
+
\sum_{f \in F}
\lambda_f f [f]
\in
\lb Y_m \rb_\ZT. 
\end{eqnarray*}
In other words $\alpha - \alpha_1 \in \lb Y_m \rb_\ZT$ where:
\begin{equation}\label{lab_doublestar}
\alpha_1 = 
\sum_{x \in X_{m-1}}
\lambda_x [x]
+
\sum_{x \in X_{m-2}}
\lambda_x [x],
\end{equation}
and
\begin{equation}
\label{eqn_TiffH2}
\sum_{x \in X_{m-1}}
\lambda_x x
+
\sum_{x \in X_{m-2}}
\lambda_x (1-e) [x] 
=
\partial_m'(\alpha_1)
=
0.
\end{equation}
Now
\[
\alpha_1
=
\left(
\sum_{x \in X_{m-1}}
\lambda_x (1-e) [x]
+
\sum_{x \in X_{m-2}}
\lambda_x e [x]
\right)
+
\alpha_2
\]
where
\[
\alpha_2 = \sum_{x \in X_{m-1}}
\mu_x [x]
+
\sum_{x \in X_{m-2}}
\sigma_x[x],
\]
with $\mu_x = \lambda_x e \in \ZL$ $(x \in X_{m-1})$, and $\sigma_x = \lambda_x (1-e) \in \mathbb{Z}$ $(x \in X_{m-2})$ since $l(1-e) = 0$ for all $l \in L$. Since 
\[
(1-e)[x] \ (x \in X_{m-1}), 
\quad 
e[x] \ (x \in X_{m-2}),
\]
all belong to $Y_m$ we deduce that $\alpha_1 - \alpha_2 \in \lb Y_m \rb_\ZT$. Applying $\partial_m'$ to $\alpha_2$ gives:
\begin{equation}
\label{eqn_TiffH3}
\sum_{x \in X_{m-1}}
\mu_x x
+
\sum_{x \in X_{m-2}}
\sigma_x (1-e)[x] = \partial_m'(\alpha_2) = 0,
\end{equation} 
where $\mu_x \in \ZL$ $(x \in X_{m-1})$ and $\sigma_x \in \mathbb{Z}$ $(x \in X_{m-2})$. 
For each $x \in X_{m-1}$, since $\mu_x \in \ZL$ it decomposes uniquely in the following way
\begin{equation}
\label{eqn_triangle_triangle}
\mu_x = \sum_{f \in F} \mu_x^{(f)}, \quad (\mu_x^{(f)} \in \ZH_f).
\end{equation}
In the next step of the argument our aim is to deduce $\sigma_x = 0$ for all $x \in X_{m-2}$. 
First we consider what happens when $m \geq 3$ (recall that $m$ is odd by assumption), and then we shall explain how to modify the argument in the case $m=1$. 

Suppose $m \geq 3$. Then, for each $x \in X_{m-2}$ by considering the coefficient of $[x]$ in equation \eqref{eqn_TiffH3}, and recalling \eqref{eqn_newstar} and the facts $\mu_x \in \ZL$ $(x \in X_{m-1})$ and $\sigma_x \in \mathbb{Z}$ $(x \in X_{m-2})$, we immediately deduce that $\sigma_x = 0$ for all $x \in X_{m-2}$.

Turning our attention to the special case when $m=1$. In this case the above argument leads to: 
\[
\alpha_2 = \sum_{x \in X_0} \mu_x [x]
+
\sum_{f \in F} \sigma_f [f]
\]
where $\mu_x \in \ZL$ $(x \in X_0)$ and $\sigma_f \in \mathbb{Z}$ $(f \in F)$, and 
\begin{equation}
\label{eqn_TiffH4}
\sum_{x \in X_0} \mu_x x
+
\sum_{f \in F} \sigma_f (1-f) 
= 
\partial_1'(\alpha_2)
=
0.
\end{equation}
For every $f \in F$ and $x \in X_0 \subseteq \ZH$ since $H$ stabilises $H_f$ under its action by right multiplication it follows that, in the notation of \eqref{eqn_triangle_triangle}, $\mu_x^{(f)} x \in \ZH_f$. Thus from equation \eqref{eqn_TiffH4} we conclude that for each $f \in f$:
\begin{equation}
\label{eqn_TiffH5}
\sum_{x \in X_0} \mu_x^{(f)} x - \sigma_f f = 0,
\end{equation}
where $\mu_x^{(f)} \in \ZH_f$ and $\sigma_f \in \mathbb{Z}$. 
But $X_0 \subseteq \ker \epsilon_T$ where $\epsilon_T$ is the standard augmentation, and so from equation \eqref{eqn_TiffH5} it follows that 
$\sigma_f f \in \ker \epsilon_T$ and so since $\sigma_f \in \mathbb{Z}$ we have
$
\sigma_f = \epsilon_T(\sigma_f f) = 0
$
for all $f \in F$. 

Thus both in the case $m = 1$ and $m \geq 3$ we conclude that
\begin{equation}
\label{**}
\alpha_2
=
\sum_{x \in X_{m-1}} \mu_x [x],
\end{equation}
where $\mu_x \in \ZL$ $(x \in X_{m-1})$, and
\[
\sum_{x \in X_{m-1}} \mu_x x
=
\partial_m'(\alpha_2)
=
0.
\]
For all $l \in L$ and $x \in X_{m-1}$, from \eqref{eqn_newstar} we deduce
\[
lx \in \bigoplus_{x \in X_{m-2}} \ZH_f [x] 
\Leftrightarrow
l \in H_f, 
\]
and it follows that
\[
\alpha_2 = \sum_{f \in F} \alpha_2^{(f)}
\]
where, using the decomposition \eqref{eqn_triangle_triangle}, $\alpha_2^{(f)} = 
\sum_{x \in X_{m-1}}
\mu_x^{(f)} [x]  \in \ker \partial_m'$ for every $f \in F$. 
Let $f \in F$ be arbitrary. Then
\[
\alpha_2^{(f)} = 
\sum_{x \in X_{m-1}}
\mu_x^{(f)} [x] 
\quad
(\mu_x^{(f)} \in \ZH_f),
\]
and
\begin{equation}
\label{eqn_triangle}
\sum_{x \in X_{m-1}} \mu_x^{(f)} x = \partial_m'(\alpha_2^{(f)}) = 0.
\end{equation}
To complete the proof of the claim it suffices to show that $\alpha_2^{(f)} \in \lb Y_m \rb_\ZT$. 
Clearly 
\[
f \cdot X_m  = f X_m \subseteq \lb Y_m \rb_\ZT.
\]
From equation \eqref{eqn_triangle} we deduce that
$
e \alpha_2^{(f)}  
$
belongs to
$
\ker \partial_m = \lb X_m \rb_\ZH.
$
So we can write:
\[
e \alpha_2^{(f)} = 
\mu_1 z_1 + \cdots + \mu_r z_r, \quad (\mu_i \in \ZH, z_i \in X_m).
\]
But then since $fe = f$ and $f$ is a left identity for $\alpha_2^{(f)}$ we obtain 
\begin{eqnarray*}
\alpha_2^{(f)}  = 
fe \alpha_2^{(f)} 
& = &
f(\mu_1 z_1 + \cdots + \mu_r z_r) \\
& = & 
(f \mu_1) (f z_1) + \cdots + (f \mu_r) (f z_r) 
 \in 
\lb f X_m \rb_\ZT \subseteq \lb Y_m \rb_\ZT, 
\end{eqnarray*}
since $f \mu_i \in \ZH_f$ $(i=1, \ldots, r)$ and $f$ is a right identity for $L_f$. 
Since $f$ was arbitrary this shows that $\alpha_2 \in \lb Y_m \rb_\ZT$, completing the proof of the claim. 
\end{proof}
Returning to the proof of Proposition~\ref{prop_diamond}, by
definition each of $B_0, B_1, \ldots, B_n$ is a free left $\ZT$-module. To complete the proof we must show that $\mathcal{B}$ is exact. Clearly $\epsilon_T$ maps $B_0$ onto $\mathbb{Z}$. 
Now let $1 \leq j \leq n-1$ and consider $\im \partial_{j+1}'$. It follows from the claim that $Y_j \subseteq \ker \partial_j'$ which, along with the definitions of $Y_j$, $B_j$ and $\partial_j'$, shows that $\im \partial_{j+1}' \subseteq \lb Y_j \rb_\ZT \subseteq \ker \partial_j'$. Since $Y_j$ is a subset of $\im \partial_{j+1}'$ it then follows from the claim that the sequence $\mathcal{B}$ is exact. 

The last clause in the statement of the proposition follows since if $X_0, \ldots, X_{m-1}$ (for some $1 \leq m \leq n$) are all finite, and $F$ is finite, then from its definition $B_m$ will clearly be finitely generated.
\end{proof}

\begin{proof}[Proof of Theorem~\ref{thm_completelysimplemain}]
Let $L$ be the $\gl$-class of $U$ that contains $H$. Then: 
\begin{align*}
&\;\phantom{\iff} 		\;\mbox{$U$ is of type left-$\FPn$} & \\
&\iff 				\mbox{$L$ is of type left-$\FPn$} & \mbox{(by Proposition~\ref{prop_SiffT})}  \\
&\iff 				\mbox{$E(L) < \infty$ \& $H$ is of type $\FPn$} & \mbox{(by Proposition~\ref{prop_TiffH})}  \\
&\iff 				\mbox{$U$ has finitely many right ideals}  \\
&\;\phantom{\iff}		\;\mbox{\& $H$ is of type $\FPn$.}
\end{align*}
\end{proof}

\section{Kobayashi's Criterion and the Property ${\rm FP}\sb 1$}
\label{sec_FP1}

In this section we turn our attention to the particular case $n=1$ and examine the behaviour of the property $\FP_1$ in more detail. As mentioned in the introduction, a group is of type $\FPone$ if and only if it is finitely generated. For monoids this is no longer the case. In a recent paper Kobayashi characterised the property left-$\FPone$ for monoids in the following way.

Let $S$ be a semigroup and let $A$ be a subset of $S$. 
A subsemigroup $T$ of $S$ is called \emph{right unitary} if $s t \in T$ implies $s \in T$ for any $t \in T$ and $s \in S$.
The intersection of two right unitary subsemigroups is clearly right unitary, so we may speak of the right unitary subsemigroup of $S$ generated by $A$, which we denote by  $\lb A \rb_{r.u.}$. We say that $S$ is \emph{right unitarily finitely generated} if there is a finite subset $A$ of $S$ such that $\lb A \rb_{r.u.} = S$.   
For a subset $A$ of a semigroup $S$ we use $\lb A \rb$ to denote the subsemigroup of $S$ generated by $A$. 

The \emph{right Cayley graph} $\Gamma_r(S,A)$ of $S$ with respect to a subset $A$ of $S$ is the directed labelled graph with vertices the elements of $S$, and a directed edge from $x$ to $y$ labelled by $a \in A$ if and only if $xa = y$ in $S$.
(Note here that we do not insist that $A$ is a generating set for $S$.) 
We write this as $x \xrightarrow{a} y$. We say that $\Gamma_r(S,A)$ is \emph{connected} if between any two vertices there is an undirected path. Kobayashi's characterisation of the property left-$\FPone$ for monoids may be stated as follows.

\begin{theorem}\cite[Proposition~2.4 \& Theorem~2.6]{Kobayashi2009}
\label{thm_Kobayashi}
A monoid $S$ is of type left-$\FP_1$ if and only if there is a finite subset $A$ of $S$ such that one of the following equivalent conditions is satisfied:
	\begin{enumerate}
	\item[(i)] $S$ is right unitarily generated by $A$ i.e. $\lb A \rb_{r.u.} = S$;
	\item[(ii)] the right Cayley graph $\Gamma_r(S,A)$ is connected. 
	\end{enumerate}
\end{theorem}

Kobayashi's characterisation allows for a more detailed analysis of $\FPn$ for monoids in the special case that $n=1$. In this section among other things we use his criterion as a tool to characterise completely simple semigroups of type left-$\FPone$ in terms of the number of right ideals and the subsemigroup generated by the idempotents. 
The following concept will play an important role. 

\begin{definition}[Relative rank]
Let $S$ be a semigroup and let $T$ be a subsemigroup of $S$. Then we define:
\[
\mathrm{rank}(S:T) = \inf_{A \subseteq S}\{ |A| : \lb T \cup A \rb = S   \}
\]
which we call the \emph{relative rank of $T$ in $S$}.
\end{definition}

It is not so surprising that the idea of relative rank arises here, since it is a notion that is central for understanding generating sets of completely simple semigroups; see \cite{Gray2005, Ruskuc1994}.

\begin{theorem}
\label{thm_fp1compsimple}
Let $U$ be a completely simple semigroup, let $H$ be a maximal subgroup of $U$ and let $\lb E(U) \rb$ be the subsemigroup generated by the idempotents of $U$. Then $U$ is of type left-$\FP_1$ if and only if the following two conditions are satisfied:
\begin{enumerate}
\item[(i)] $U$ has finitely many right ideals, and
\item[(ii)] the subgroup $K$ of $H$ generated by $\lb E(S) \rb \cap H$ has finite relative rank in $H$. 
\end{enumerate} 
\end{theorem}
\begin{proof}
Let $S = U^1$. 
By the Rees theorem (see \cite[Chapter~3]{Howie1995}, or originally \cite{Rees1940}) we may identify $U$ with a Rees matrix semigroup $M[G;I,\Omega;P]$ over a group $G$ where $G \cong H$ and $P = (p_{\omega i})_{\omega \in \Omega, i \in I}$ is an $\Omega \times I$ matrix with entries from $G$.  This semigroup has elements $U = I \times G \times \Omega$ and multiplication given by:
\[
(i,g,\omega) (j, h, \mu) = (i, g p_{\omega j} h, \mu).
\]
Moreover, we may assume that the matrix $P$ is in normal form i.e. $p_{1 \omega} = p_{i 1} = 1$ for all $i \in I$ and $\omega \in \Omega$ (see \cite[Chapter~3, Section~4]{Howie1995}) and we may also suppose without loss of generality that $H = H_{1 1}$. It is well known and easy to prove (see \cite{Howie1978} for example) that the subgroup
$K$ of $H$ generated by $\lb E(S) \rb \cap H$ consists of all triples $(1,g,1)$ where $g$ belongs to the subgroup of $G$ generated by the entries in the matrix $P$.

Suppose that $I$ is finite (i.e. that $U$ has finitely many right ideals) and $K$ has finite relative rank in $H$. Let $X$ be a finite subset of $H$ such that $\lb X \cup K \rb = H$. Let $F$ be the set of all idempotents in some fixed $\gl$-class $L$ of $U$. We claim that $F \cup X$ is a right unitary generating set for $U$. First observe that $E(U) \subseteq \lb F\cup X \rb_{r.u}$. Indeed, given any $e \in E(U)$ there exists $f \in F$ such that $e \gr f$ and so $ef = f$ and thus $e \in \langle F \cup X \rangle_{r,u}$. Therefore
\begin{eqnarray*}
\lb F \cup X \rb_{r.u}
& = &
\lb F \cup X \cup E(U) \rb_{r.u}
 =
\lb X \cup E(U) \rb_{r.u} \\
& \supseteq &
\lb X \cup E(U) \rb
 \supseteq
\lb (X \cup (\lb E(U) \rb  \cap H)) \cup E(U) \rb \\
& = &
\lb H \cup E(U) \rb = U.
\end{eqnarray*}
Hence $\lb F \cup X \rb_{r.u} = S$ and since $F \cup X$ is finite it follows from Theorem~\ref{thm_Kobayashi} that $S$ is of type left-$\FPone$.

For the converse, suppose that $S$ is of type left-${\rm FP}\sb 1$. Let $A$ be a finite subset of $S$ such that $\Gamma_r(S,A)$ is connected. Also, we may assume that $A$ is chosen so that 
$(i,g,\omega) \in A$ if and only if $(i, g^{-1}, \omega) \in A$ for all $i \in I$, $\omega \in \Omega$ and $g \in G$. Since for any collection of $\gr$-classes $R_1, \ldots, R_k$ of $U$ the union $R_1 \cup \ldots \cup R_k \cup \{ 1 \}$ is a right unitary submonoid of $S$ it follows that $A$ must intersect every $\gr$-class of $S$, and thus $I$ is finite.
To complete the proof we have to show that $K$ has finite relative rank in $H$.

We claim that $\lb A \cup E(S) \rb = S$. To see this, suppose for the sake of a contradiction that
$\lb A \cup E(S) \rb \subsetneq S$. Since $\Gamma_r(S,A)$ is connected (as an undirected graph)
it follows that there exist $u,v \in S$ such that $u \in \lb A \cup E(S) \rb$, $v \not\in \lb A \cup E(S) \rb$ and $u$ adjacent to $v$ in $\Gamma_r(S,A)$. So there exists $a \in A$ such that either $ua = v$ or $va = u$.
Since $v \not\in \lb A \cup E(S) \rb$ we cannot have $ua = v$ and therefore must have $va  = u$. 
Clearly none of $a$, $v$ or $u$ is equal to $1$. 
Let
\[
v = (i,g,\omega), \quad a = (j, \alpha, \mu), \quad u = (k,\beta, \nu),
\]
where $i,j,k \in I$, $\omega, \mu, \nu \in \Omega$ and $g, \alpha, \beta \in G$. 
Then
\[
(i,gp_{\omega j} \alpha, \mu) = (i,g,\omega)(j, \alpha, \mu) = va = u = (k,\beta, \nu),
\]
and so
\[
i=k, \quad \mu = \nu, \quad \& \quad g = \beta \alpha^{-1} p_{\omega j}^{-1}.
\]
Therefore since the matrix $P$ is in normal form:
\begin{eqnarray*}
v=
(i,g,\omega)
& = &
(i, \beta \alpha^{-1} p_{\omega j}^{-1}, \omega) 
 = 
(i, \beta p_{\mu 1} 1 p_{1 j } \alpha^{-1} p_{\mu 1} 1 p_{1 j} p_{\omega j}^{-1}, \omega) \\
& = &
(i, \beta, \mu)(1,1,1)(j,\alpha^{-1},\mu)(1,1,1)(j,p_{\omega j}^{-1}, \omega)
\in
\lb A \cup E(S) \rb,
\end{eqnarray*}
since $(i,\beta,\mu) = (k, \beta, \nu) = u \in \lb A \cup E(S) \rb$, $(j, \alpha^{-1}, \mu) \in A$
and $(j,p_{\omega j}^{-1}, \omega) \in E(S)$. This
is a contradiction. We conclude that $\lb A \cup E(S) \rb = S$. In particular
$H \subseteq \lb A \cup E(S) \rb$ so given $g \in G$ we can write
\[
(1,g,1) = (i_1,g_1,\omega_1)(i_2,g_2,\omega_2) \cdots (i_k,g_k,\omega_k)
\]
where $i_1 = 1$, $\omega_k = 1$ and
each $(i_r, g_r, \omega_r) \in A \cup E(S)$. It follows that
\[
g = g_1 p_{\omega_1 i_2} g_2 \cdots p_{\omega_{k-1} i_k} g_k
\]
where each $g_i$ is either the inverse in $G$ of an entry from $P$, or is the middle entry of some triple from $A$. It follows that if $B$ is the set of all middle entries of elements of $A$ then $N \cup B$ generates $G$ where $N$ is the subgroup of $G$ generated by the entries in the matrix $P$. Since $A$ is finite it follows that $B$ is finite, therefore $K$ has finite relative rank in $H$ and this completes the proof of the theorem.
\end{proof}

We leave as an open question the problem of extending Theorem~\ref{thm_fp1compsimple} to values of $n$ greater than $1$.
It seems likely that the formulation of such a result will need the introduction of the notion of a subgroup $K$ being of \emph{relative type-$\FP_n$} in a group $G$.

The following result shows how left-$\FPone$ holding in an ideal of a semigroup influences the same property holding in the semigroup.

\begin{proposition}
\label{prop_ideal}
Let $S$ be a monoid and let $J$ be a left ideal of $S$. If $J$ is of type left-$\FPone$ then so is $S$.
\end{proposition}
\begin{proof}
Let $T = J \cup \{1 \} \leq S$. 
Since $J$ is left-$\FPone$ there is a non-empty finite subset $A$ of $T$ such that $\Gamma_r(T,A)$ is connected. 
Let $a \in A \cap J$. Then for all $s \in S$ we have $sa \in J$, since $J$ is a left ideal, and it follows that $\Gamma_r(S,A)$ is connected and hence $S$ is of type left-$\FPone$.
\end{proof}

The following example shows that the converse of this result fails, even for two-sided ideals and as a consequence also shows that Theorem~\ref{thm_IdealWithIdentity} above does not hold if we remove the assumption that the ideal has a two-sided identity element.
Moreover, it shows that the finiteness assumption on $J$ in Theorem~\ref{thm_fp1minimalideal} below really is necessary.

\begin{example}
\label{ex_IdealsInGeneral}
Let $S$ be a finitely generated monoid with a minimal ideal $R$ where $R$ is isomorphic to a rectangular band $A \times B$ where $A$ and $B$ are both infinite sets. Since $S$ is finitely generated it is, in particular, of type left-$\FPone$. However, by Theorem~\ref{thm_Kobayashi}, $R$ is neither of type left-$\FPone$ nor right-$\FPone$.

Such an example $S$ may be constructed in the following way. First let $U$ be the submonoid of the full transformation monoid $\mathcal{T}(\mathbb{N})$ (the semigroup of all maps from $\mathbb{N}$ to $\mathbb{N}$ under composition) generated by $\alpha$ which maps $i \mapsto i+1$ $(i \in \mathbb{N})$ and the constant map $\gamma$ with image  $1$. The monoid $U$ is finitely generated, and has a minimal ideal isomorphic to a right zero semigroup. Let $V$ be the dual of $U$, which is finitely generated and has a minimal ideal that is an infinite left zero semigroup. Then define $S = U \times V$ which has a minimal ideal that is a rectangular band with infinitely many left and right ideals. Also $S$ is finitely generated since it is a direct product of two finitely generated monoids.
\end{example}

Combining Theorem~\ref{thm_fp1compsimple} and Proposition~\ref{prop_ideal} we obtain the following.

\begin{theorem}
\label{thm_fp1minimalideal}
Let $S$ be a monoid with a minimal ideal $J$ that is completely simple and has finitely many right ideals.
Let $H$ be a maximal subgroup of $S$ in $J$ and $K$ be the subgroup of $H$ generated by $\lb E(J) \rb \cap H$.
Then $S$ is of type left-$\FPone$ if and only if $K$ has finite relative rank in $H$.
\end{theorem}
\begin{proof}
One direction is a direct consequence of Theorem~\ref{thm_fp1compsimple} and Proposition~\ref{prop_ideal}.

For the converse, suppose that $S$ is of type left-$\FPone$. By Theorem~\ref{thm_fp1compsimple} to complete the proof it suffices to show that $T = J \cup \{ 1 \} \leq S$ is of type left-$\FPone$. Let $A$ be a finite subset of $S$ such that 
$\Gamma_r(S,A)$ is connected as an undirected graph.
Fix an $\gl$-class $L$ of $J$ and let $F = E(S) \cap L$ which is finite by assumption.
Define $B = F \cup AF$ which is a finite subset of $J$ since $J$ is an ideal.
We shall now prove that $B$ is a right unitary generating set for $T = J \cup \{ 1 \}$.

Let $e \in F$ be arbitrary and let $R = R_e$ be its $\gr$-class. We claim that for every $x \in R$ there is a path in $\Gamma_r(T,B)$ from $x$ to $e$.
Let $x \in R$. 
Since $\Gamma_r(S,A)$ is connected it follows that there is a sequence
\[
1 = y_0, y_1, \ldots, y_r = x
\]
of elements of $S$ such that for all $i$, $y_i$ and $y_{i+1}$ are connected by an arc (in some direction) in $\Gamma_r(S,A)$. Now consider the sequence:
\[
z_0 = e = e1 = ey_0, \; z_1 = ey_1, \ldots, \; z_r = ey_r = ex = x,
\]
recalling that $e$ is a left identity in its $\gr$-class. Since $\gr$ is minimal we have $ey_i \in R$ for all $i$. Now for any $a \in A$, $y_i a = y_{i+1}$ implies $(ey_i) a = (e y_{i+1})$, while $y_{i+1} a = y_i$ implies $(e y_{i+1}) a = (e y_i)$. So the sequence $(z_i)_{0 \leq i \leq r}$ is a path in $\Gamma(S,A)$ contained in $R$, beginning at $e$ and terminating at $x$. Consider a typical arc in this path: $ua = v$. This implies $u(ae) = v(e)$ and hence, since $ae, e \in B$ by definition, $u$ and $v$ are joined by a path of length at most $2$ in $\Gamma_r(T,B)$. We conclude that in $\Gamma_r(T,B)$ there is a path from every vertex $x \in R$ to $e$. Since $R$ was an arbitrary $\gr$-class, the same is true for every $\gr$-class of $J$. Also every pair of idempotents $e_i$ and $e_j$ of $F$ are connected
by a path in $\Gamma_r(T,B)$ of length $2$ via $1 \in T$ since $e_i, e_j \in B$.
We conclude that $\Gamma_r(T,B)$ is connected and so $T$ is of type left-${\rm FP}\sb 1$.
\end{proof}

In particular, in Theorem~\ref{thm_fp1minimalideal}, if $E(J)$ is finite then $K$ has finite relative rank in $H$ if and only if $H$ is finitely generated which gives the following result. 

\begin{corollary}
\label{cor_fmanyLandRclasses}
Let $S$ be a semigroup with finitely many left and right ideals, let $K$ be the unique minimal ideal of $S$, and let $H$ be a maximal subgroup of $S$ in $K$. Then the following are equivalent:
\begin{enumerate}
\item[(i)] $S$ is of type left-$\FP_1$;
\item[(ii)] $S$ is of type right-$\FP_1$;
\item[(iii)] the group $H$ is finitely generated (equivalently, $H$ is of type $\FP_1$). 
\end{enumerate}
\end{corollary}
Currently we do not know whether Corollary~\ref{cor_fmanyLandRclasses} holds for left-$\FPn$ for values of $n$ greater than one. 
We do however know that it holds in one direction, passing from $S$ to $H$, by virtue of Theorem~\ref{thm_MinimalMaximalSubgroups}.  

\section{Further Applications and Examples}
\label{sec_applications}

In this section we give some examples showing that the finiteness conditions imposed in the main results of the paper really are necessary. We also give some further applications of our results. 

The following example shows that without the finiteness assumption on the number of left ideals, Theorem~\ref{thm_MinimalMaximalSubgroups} (and therefore also Theorem~\ref{thm_completelysimplemain}) no longer holds.  

\begin{example}
\label{ex_NeedBFinite}
Let $G$ be the free group over $X$ where $X = \{ x_i : i \in \mathbb{N} \}$. Let $S = M[G;I,\Omega;P]$ 
be the Rees matrix semgiroup over $G$ with structure matrix
\[
P =
\begin{pmatrix}
1 & 1 \\
1 & x_1 \\
1 & x_2 \\
\vdots & \vdots
\end{pmatrix}.
\]
Then by Theorem~\ref{thm_fp1compsimple}, $S$ is of type left-${\rm FP}\sb 1$. But $G$ is not of type $\FPone$ since $G$ is an infinitely generated group. 
\end{example}

\subsection*{Simple semigroups} In Theorem~\ref{thm_main2} we proved that if a completely simple semigroup is of type left- and right-$\FP_n$ then all of its maximal subgroups are of type $\FP_n$. We shall now see that this result does not extend to simple semigroups in general.

We begin by quoting a well-known result regarding $\FP_n$ for amalgamated free products of groups; see \cite{Bieri1976} for a proof.

\begin{proposition}
\label{prop_freeproduct}
Let $G$ be the amalgamated free product $A *_C B$ of groups $A$,$B$ and $C$, and let $n \in \mathbb{N}$.
If $G$ is $\FP_n$ and $A$ and $B$ are $\FP_{n-1}$ then $C$ is $\FP_{n-1}$.
\end{proposition}

\begin{example}
Let $S$ be the monoid defined by the following finite presentation:
\[
\begin{array}{c}
\langle
a_1, a_2, a_3, a_4, \
a_1', a_2', a_3', a_4', \
b,c  \ | \
a_j a_j' = a_j' a_j = 1, \ a_1 a_2 = a_3 a_4, \ bc=1, \\
b a_j = a_j^2 b, \ a_j c = c a_j^2 \ (j=1,2,3,4)
\rangle.
\end{array}
\]
Then from \cite[Proposition~3.3]{Ruskuc1999} the group of $G$ units of $S$ is defined by the following group presentation:
\[
\langle
a_1, a_2, a_3, a_4
\ | \
a_1^{2^i} a_2^{2^i} = a_3^{2^i} a_4^{2^i} \
(i=0,1,2,3, \ldots )
\rangle.
\]
As observed in \cite{Ruskuc1999} $S$ is a Bruck--Reilly extension (see \cite[Chapter~5]{Howie1995} for a definition of Bruck--Reilly extension) of the group $G$  and consequently $S$ is simple and every maximal subgroup of $S$ is isomorphic to $G$.
Also, clearly $G$ is the free product with amalgamation of $A_1$, $A_2$ (both free groups of rank $2$) over a non finitely generated subgroup. Therefore it follows by Proposition~\ref{prop_freeproduct} that $G$ is not of type $\FP_2$. Recalling (see \cite[Proposition~5.6]{Otto1997} for instance) that every finitely presented monoid is of type left- and right-$\FP_2$ we obtain the following. 

\begin{proposition}
There exists a simple monoid that is finitely presented, and hence of type left- and right-$\FP_2$, none of whose maximal subgroups are of type $\FP_2$.
\end{proposition}

This leads naturally to the following question: is it true that for every $n \geq 1$ there is a simple monoid of type left- and right-$\FPn$ none of whose maximal subgroups are of type $\FPn$? 

\end{example}

\subsection*{Strong semilattice of monoids}

We may apply our results to another fundamental construction from semigroup theory, so-called \emph{strong semilattices of monoids}.
Let $Y = (Y, \leq)$ be a semilattice, and let $A_{\alpha} (\alpha \in Y)$, be a family of disjoint monoids indexed by $Y$. Denote by $1_{\alpha}$ the identity of $A_{\alpha}$. Suppose that for any two elements $\alpha, \beta \in Y$, $\beta \leq \alpha$, there exists a homomorphism $\phi_{\alpha, \beta} : A_{\alpha} \rightarrow A_{\beta}$ such that
\begin{enumerate}
\item $\phi_{\alpha, \alpha}$ is the identity homomorphism on $A_{\alpha}$
\item $\phi_{\alpha, \beta} \phi_{\beta, \gamma} = \phi_{\alpha,\gamma}$, for all $\alpha, \beta, \gamma \in Y$ with $\gamma \leq \beta  \leq  \alpha$.
\end{enumerate}
The set $S = \cup_{\alpha \in Y} A_{\alpha}$ can then be made into a semigroup by defining
\[
ab = (a \phi_{\alpha, \alpha \beta})(b \phi_{\beta, \alpha \beta}), \quad a \in A_{\alpha}, \; b \in A_{\beta}.
\]
When all $A_{\alpha}$ are groups then we obtain exactly the Clifford monoids (originally introduced in \cite{Clifford1941}) from this construction (a \emph{Clifford monoid} is a regular semigroup whose idempotents are central). More details on this construction may be found in \cite[Chapter~4]{Howie1995}.

The following result characterises the property left-$\FP_n$ for strong semilattices of monoids, and so in particular for Clifford monoids. 

\begin{theorem}
\label{thm_Clifford}
Let $S = \mathcal{S}[Y;A_{\alpha},\phi_{\alpha,\beta}]$ be a strong semilattice of monoids. Then $S$ is of type left-${\rm FP}\sb n$ if and only if $Y$ has a minimal element $e$, and the monoid $A_e$ is of type left-${\rm FP}\sb n$.
\end{theorem}
\begin{proof}
Suppose $S$ is of type left-${\rm FP}\sb n$. Then in particular $S$ is of type left-${\rm FP}\sb 1$ and hence the semilattice $Y$, which is a retract of $S$, is also of type ${\rm FP}\sb 1$ by \cite[Theorem~3]{Pride2006}. (In fact, left-${\rm FP}\sb 1$ is even preserved by arbitrary homomorphic images, which is easily seen from Theorem~\ref{thm_Kobayashi}.) By Theorem~\ref{thm_Kobayashi} this means that $Y$ is right unitarily finitely generated. Let $A$ be a finite right unitary generating set for $Y$.
Let $X$ be the subsemilattice generated by $A$, which is finite since $A$ is finite, and let $z$ be the unique minimal element of $X$.
Define $Z = \{ y \in Y : zy = yz = z  \}$. Clearly $Z$ is a subsemigroup of $S$, and is right unitary since for $x \in Z$ and $y \in Y$, $yx \in Z$ implies $yxz = z$ and so $yz = z$ which gives $y \in Z$. Since $A$ is a right unitary generating set for $Y$ it follows that $Z = Y$ and so $z$ is the unique minimal element of $Y$. 
The result now follows by applying Theorem~\ref{thm_IdealWithIdentity}.
\end{proof}

For inverse semigroups it is known that left-${\rm FP}\sb n$ and right-${\rm FP}\sb n$ are equivalent. In particular this is true for Clifford monoids. The above theorem applies to Clifford monoids.

\begin{corollary}
A Clifford monoid is of type $\FPn$ if and only if it has a minimal idempotent $e$ and the maximal subgroup $G_e$ containing $e$ is of type $\FPn$. 
\end{corollary} 

Another related application of the results of Section~\ref{sec_submonoids} is the following. 

\begin{corollary}
Let $S$ be an inverse semigroup with a minimal idempotent $e$, and let $G$ be the maximal subgroup of $S$ containing $e$. Then $S$ is of type $\FPn$ if and only if $G$ is of type $\FPn$.
\end{corollary}

\section{Other homological finiteness properties}
\label{sec_other_props}

We conclude the paper with some remarks about some other homological finiteness properties of monoids. 

\subsection*{\boldmath The property bi-$\FP_n$.}
In \cite{Alonso2003} Alonso and Hermiller introduced a property which they called \emph{bi-$\FP_n$} (the same property is called \emph{weak bi-$\FP_n$} in \cite{Pride2006}).

A monoid $M$ is said to be of type \emph{bi-$\FP_n$} if there is a finite rank length $n$ resolution of $\mathbb{Z}$ by $(\ZM,\ZM)$-bimodules. 

Pride \cite{Pride2006} showed that a monoid is of type bi-$\FP_n$ (in the sense of Alonso and Hermiller) if and only if it is of type left- and right-$\FP_n$. Therefore an alternative way or expressing Theorem~\ref{thm_main2} above is as follows. 

\begin{theorem}
A completely simple semigroup $S$ is of type bi-$\FP_n$ if and only if it has finitely many left and right ideals and all of its maximal subgroups are of type $\FP_n$.
\end{theorem}

\subsection*{Cohomological dimension}

Several problems regarding closure properties of homological finiteness conditions of monoids were posed in \cite[Remark and Open Problem 4.5]{Ruskuc1999} and in \cite[Open Problem 11.1(i)]{Ruskuc1998}.
Specifically, in \cite[Remark and Open Problem 4.5]{Ruskuc1999} it was asked whether for a regular semigroup $S$ with finitely many left and right ideals whether $S$ has property left-$\FPn$ (resp. finite co-homological dimension) if and only if all maximal subgroups of $S$ have property left-$\FPn$ (resp. finite co-homological dimension). We have already observed above that the first of these questions, concerning property left-$\FPn$, has a negative answer. We may similarly answer negatively the question about cohomological dimension using \cite[Theorem~1]{Guba1998} which states that a monoid has left and right cohomological dimension zero if and only if it has a two-sided zero element. Therefore, by taking a group $G$ with infinite cohomological dimension and adjoining a zero element we obtain a regular monoid with finitely many left and right ideals, and with finite cohomological dimension, but with a maximal subgroup that has infinite cohomological dimension.  

In exactly the same way we see that neither the property left-$\FPn$, nor that of having finite cohomoligical dimension, is inherited by subsemigroups with finite Rees index (the Rees index of a subsemigroup $T$ of a semigroup $S$ is defined as the cardinality of the complement $S \setminus T$). This answers negatively two further open problems that were posed in \cite[Open Problem 11.1(i)]{Ruskuc1998}. 

\section*{Acknowledgements}

The authors would like to thank Professor Benjamin Steinberg for useful discussions during the preparation of this paper.

\bibliographystyle{abbrv}
\def\cprime{$'$} \def\cprime{$'$} \def\cprime{$'$}

\end{document}